\documentclass{article}
\usepackage{graphicx} 
\usepackage{amsmath}
\usepackage{amssymb}
\usepackage{amsthm}
\usepackage{amsfonts}
\usepackage{array}
\usepackage[x11names,dvipsnames]{xcolor}
\usepackage{pgf,tikz}

\usepackage{pgfplots}
\usepgfplotslibrary{fillbetween}

\usetikzlibrary{arrows}
\usetikzlibrary{hobby}

\usepackage{tikz,tkz-euclide}

\usepackage{url}
\usepackage{hyperref}

\newtheorem{thmX}{Theorem}

\newtheorem{theorem}{Theorem}
\newtheorem*{thmJ*}{Theorem J'}
\newtheorem{proposition}{Proposition}[section]
\newtheorem{definition}[proposition]{Definition}

\newtheorem{lemma}[proposition]{Lemma}

\newtheorem{corollary}[proposition]{Corollary}
\newtheorem{remark}[proposition]{Remark}
\newtheorem{example}[proposition]{Example}

\newcommand{\black}{\color{black}{}}

\def\bN{{\mathbb N}}
\def\bR{{\mathbb R}}
\def\bZ{{\mathbb Z}}
\def\cA{{\cal A}}

\def\cC{{\cal C}}

\def\cF{{\cal F}}

\def\cL{{\cal L}}
\def\fL{{\mathfrak L}}
\def\cO{{\cal O}}
\def\NW{{\cal P}}
\def\cR{{\cal R}}
\def\cS{{\cal S}}

\def\Down{\text{Down}}

\newcommand{\BF}{\bf\boldmath }
\def\NWto{\succcurlyeq_{{\cal P}_F}}
\def\Ato{\succcurlyeq_{{\cal A}_F}}

\def\NWGto{\succcurlyeq_{{\cal P}_{G}}}

\def\eps{\varepsilon}
\def\NWto{\succcurlyeq_{{\cal P}_F}}
\def\prec{\succ}

\def\vectorflow{regular\ }
\def\Rt{\Bbb R^2}

\title{Wandering Flows on the Plane}
\author{Joseph Auslander\\ Roberto De Leo }
\date{\today}

\begin{document}

\maketitle

\begin{abstract}    
    We study planar flows without non-wandering points and prove several properties of these flows in relation with their prolongational relation.
    The main results of this article are that a planar (regular) wandering flow has no generalized recurrence and has only two topological invariants: the space of its orbits and its prolongational relation (or, equivalently, its smallest stream).
    As a byproduct, our results show that, even in absence of any type of recurrence, the stream of a flow contains fundamental information on its behavior.
\end{abstract}

\section{Introduction}

Very little is available in literature on the qualitative dynamics of flows without non-wandering points.
We call such flows {\em wandering}.
Since in a compact space every orbit has limit points, that are necessarily non-wandering, wandering flows can arise only on non-compact spaces.
In this article we study the simplest non-trivial case of wandering flows, namely wandering flows on the plane.
Notice that these flows have no periodic orbits and no fixed points: each orbit goes to infinity for both $t\to\pm\infty$.
In this article, we restrict our study to {\em regular} flows, namely flows that are locally topologically equivalent to the flow of a constant vector field.

This problem has been already addressed in literature in different but equivalent forms.
In~\cite{Kap40,Kap41}, Kaplan  studied and classified regular families of curves filling the plane, namely families of planar disjoint curves that are locally homeomorphic to the family of parallel straight lines and whose union is the whole plane. 
%
%
There are two natural ways to generate a regular family of curves filling the plane: one is as the family of level sets of a continuous function without local extrema (such functions are called {\em pseudoharmonic} in literature); the other is as the family of orbits of a regular planar flow.
Kaplan, based on previous results by Whitney, showed that these three kinds of families are actually one and the same:
%
%
%
\begin{thmX}[Kaplan, 1940~\cite{Kap40}]
    Each regular family of curves filling the plane is the family of level sets of a pseudoharmonic function.
\end{thmX}
%
%
%
\begin{thmX}[Kaplan, 1948~\cite{Kap48}]
    For each 
    regular family of curves filling the plane $R$
    there is a homeomorphism of the plane onto the interior on the unit disc $D$ such that $R$
    is transformed in the family of level sets of a harmonic function on $D$.
\end{thmX}
%
%
\begin{thmX}[Whitney, 1933~\cite{Whi33}, Theorem 27A; Kaplan, 1940~\cite{Kap40}, Corollary to Theorem 42]
    \label{thm: WK}
    Every regular family of curves filling the plane is orientable and is the set of orbits of a regular planar flow.
\end{thmX}
We say that a flow has a {\em first-integral} when there is a pseudoharmonic function $H$ such that each orbit lies in a level set of $H$, namely $H$ is constant over every orbit.
The results above imply the following important result.
%
\begin{thmX}[Whitney, Kaplan]
    \label{thm: first-integral}
    Every regular planar flow has a first-integral.
\end{thmX}

In \cite{HR57}, Haefliger and Reeb developed a theory of non-Hausdorff manifolds and used it to classify planar (regular) foliations, which are the same thing as regular family of curves filling the plane.
By Theorem~\ref{thm: WK}, this is the same as the space of orbits of a regular flow on the plane.
Their elegant approach can be used to solve our problem but here, in order to keep technicalities at a minimum,  here we rather prefer using the Kaplan approach.

\medskip
An important byproduct of our results is providing a first concrete example of a claim made in~\cite{DLY24} by the second author jointly with Jim Yorke.
In that article the authors study in some detail the ``prolongation relation'' and introduce the concept of ``stream'' of a (semi-)flow%
, which are closed and transitive relations containing the prolongational relation, and argue that these relations are a concept more primitive than recurrent points and their generalizations 
such as non-wandering points, generalized recurrent points and chain-recurrent points. 
In this article we show that, although no point is recurrent or even just generalized recurrent, and so recurrence plays no role, yet the non-wandering relation (or, equivalently, the smallest stream) allows to distinguish between topologically inequivalent flows.

\medskip
In more detail, we prove the following main results.
Let $F$ be a regular planar flow. Then:
\begin{enumerate}
     \item All prolongational limit sets of $F$ of order larger than 2 are equal to the corresponding prolongational limit sets of $F$ of order 2 (Theorem~\ref{thm1}, see also Definitions~\ref{def: L^1} and~\ref{def: L^2}).
    \item $F$ does not have generalized recurrent points and, consequently, it has a $C^0$ strict Lyapunov function (Theorem~\ref{thm: A_F is empty} and Corollary~\ref{cor: Lyap}, see also Definitions~\ref{def: gener rec pts} and~\ref{def: Lyap}).
    \item $F$ is completely characterized, modulo topological equivalence, by its following two topological invariants: the topology of its space of orbits and its prolongational relation or, equivalently, its Auslander stream (Theorem\ref{thm3} and~3', see also Definitions~\ref{def: Prol rel} and~\ref{def: stream} and Proposition~\ref{prop: A_F}).
\end{enumerate}



The article is structured as follows.
In Section~\ref{sec: tools}, we introduce the main concepts and tools needed for our purposes.
In Section~\ref{sec: wandering}, we go over general properties of regular flows and provide several examples.
Finally, in Section~\ref{sec: results}, we provide a full classification of regular planar flows.

\black
\section{Basic definitions and properties}
\label{sec: tools}
Throughout this article, we denote by $X$ a metrizable locally compact connected topological space and by
$F$ a continuous-time flow on $X$, namely a continuous map $F:\bR\times X\to X$ such that $F^t(F^s(x))=F^{t+s}(x)$ for all $t,s\in\bR$ and $F^0(x)=x$.
In this section we go over some general property and tools relative to flows.
\begin{definition}
    \label{def: Prol rel}
    Given a point $p\in X$, the {\bf backward orbit}, {\bf forward orbit} and {\bf orbit} of $p$ under $F$ are the sets
    $$
    \cO^-_F(p) = \{F^t(p):t\leq0\},
    $$
    $$
    \cO^+_F(p) = \{F^t(p):t\geq0\},
    $$
    $$
    \cO_F(p) = \cO^-_F(p)\cup\cO^+_F(p).
    $$
    We call {\bf orbit space} of $F$ the binary relation
    $$
    \cO_F = \{(p,q):q\in\cO^+_F(p)\}.
    $$
    We call {\bf prolongational relation}~\cite{AG97} of $F$, denoted by {\BF$\NW_F$}, the topological closure of $\cO_F$ in $X\times X$ endowed with the product topology.
    If $(x,y)\in\NW_F$, we write {\BF $x\NWto y$} and we set
    $$\text{\BF$\Down_{\NW_F}(p)$}=\{q:p\NWto q\}\subset X.$$
    Sometimes we write that {\BF $q$ is prolongationally downstream from $p$} to say that $p\NWto q$.
\end{definition}
From now on, we will denote by $d$ a distance on $X$ compatible with its topology.
All results in the remainder of the article will not depend on the particular choice of $d$. 
The reader can verify that $p\NWto q$ if and only if, for every $\eps>0$, there is a $x$ and $t>0$ such that $d(x,p)<\eps$ and $d(F^t(x),q)<\eps$.
\begin{example}
    \label{ex: circle}
    Given $x,y\in\Bbb S^1$ with $x\neq y$, denote by $\overset{\frown}{xy}$ the closed arc that is the closure of all points $z$ such that $x,z,y$ are in clockwise order.
    Consider the flow $F$ on $\Bbb S^1$ which us three fixed points, in clockwise order, $p,q,r$ and where every other point moves clockwise.
    Then $x\NWto y$ if and only if either $\overset{\frown}{xy}\subset\overset{\frown}{pq}$ or $\overset{\frown}{xy}\subset\overset{\frown}{qr}$ or $\overset{\frown}{xy}\subset\overset{\frown}{rp}$.
    For instance, $p\NWto q$ because one can jump from $p$ to a point $x$ arbitrarily close to $p$ and use the flow to get arbitrarily close to $q$.
    Let now $G$ be an extension of $F$ to the unit disc $D$ so that every point of $D$, except for the center of the disc, moves in a spiral that asymptotes to the unit circle.
    In this case, $x\NWGto y$ for any $x,y\in\Bbb S^1$ because from any $x$ on the unit circle one can jump to a spiral arbitrarily close to $x$ and then use the flow to move on the spiral arbitrarily close to any other point $y$ on the unit circle.
\end{example}
%
\noindent
%
%
\begin{definition}
    We call {\BF$\omega$-limit set} of $p$ the set
    $$
    \text{\BF$\omega_F(p)$}=\{q:\text{there is }t_n\to\infty\text{ s.t. }F^{t_n}(p)\to q\}.
    $$
    Analogously, we call {\BF$\alpha$-limit set} of $p$ the set
    $$
    \text{\BF$\alpha_F(p)$}=\{q:\text{there is }t_n\to-\infty\text{ s.t. }F^{t_n}(p)\to q\}.
    $$
\end{definition}
The following two results illustrate some important relation between the prolongation relation and limit sets.
\begin{lemma}[De Leo and Yorke~\cite{DLY24}]
    Let $p\NWto q$. Then $p\NWto\cO^+_F(q)\cup\omega_F(q)$.
\end{lemma}
\begin{proposition}[De Leo and Yorke~\cite{DLY24}]
    \label{prop:Down}
    For every flow $F$, $\NW_F$ is invariant under the induced action of $F$ on $X\times X$ given by $F^t(x,y)=(F^t(x),F^t(y))$.
    Moreover, for every $p$, the following holds:
    \begin{enumerate}
        \item $\Down_{\NW_F}(p)$ is forward-invariant under $F$;
        \item $\Down_{\NW_F}(p)
        \;\supset\; 
        \cO^+_F(p)\cup\omega_F(p)$;
        \item if $\omega_F(p)\neq\emptyset$, 
        $\cO^+_F(p)\cup\Down_{\NW_F}(\omega_F(p))
        \supset
        \Down_{\NW_F}(p)$;
        \item 
        $\Down_{\NW_F}(p)\;=\;\cO^+_F(p)\cup\Down_{\NW_F}(q)$ for all $q\in\cO^+_F(p)$.
    \end{enumerate}
\end{proposition}

\medskip\noindent
{\bf Prolongations and prolongational limit sets.} 
As we show below, the prolongational relation can be defined also through the theory of prolongations, developed by Ura~\cite{Ura53} and Seibert and Auslander~\cite{AS64}.
In particular, the prolongational limit set will play an important role in our study of planar wandering flows.
%
\begin{definition}
    \label{def: L^1}
    Given points $p,q\in X$ and reals $\eps,T>0$, a {\BF $(\eps,T)$-link} from $p$ to $q$ is a segment of trajectory $\{F^t(x): t\in[0,\tau]\}$ such that:
    \begin{enumerate}
        \item $d(p,x)<\eps$;
        \item $d(q,F^\tau(x))<\eps$;
        \item $\tau\geq T$.
    \end{enumerate}
    %
     
%

    The {\bf first prolongation} of $p$ under $F$ is the set {\BF $D^1_F(p)$} of all points $q$ for which, for every $\eps>0$, there is a $(\eps,T_\eps)$-link from $p$ to $q$ for some $T_\eps>0$.
    The {\bf first prolongational limit set} of $p$ under $F$ is the set {\BF $\Lambda^1_F(p)$} of all points $q\in D^1_F(p)$ for which, for every $\eps>0$ and $T>0$, there is an $(\eps,T)$-link from $p$ to $q$.
%
    %
    Given $A\subset X$, we set
    $D^1_F(A)=\cup_{x\in A}D^1_F(x)$ and $\Lambda^1_F(A)=\cup_{x\in A}\Lambda^1_F(x)$.
\end{definition}
\begin{proposition}
    \label{prop: Down}
    For all $p\in X$, the following hold:
    \begin{enumerate}
        \item $D^1_F(p)=\Down_{\NW_F}(p)$;
        \item $\Down_{\NW_F}(p)=\cO^+_F(p)\cup\Lambda^1_F(p)$.
    \end{enumerate}
\end{proposition}
\begin{proof}
    (1)
    We leave the proof to the reader.
    (2)
    Since $\NW_F=\overline{\cO_F}$, it follows immediately that $\Down_{\NW_F}(p)\supset\cO^+_F(p)$.
    Assume now that $q\in D^1_F(p)\setminus\cO^+_F(p)$.
    Then there is a sequence $(p_n,q_n)\to(p,q)$ such that $q_n\in\cO^+_F(p_n)$, namely $q_n=F^{t_n}(p_n)$ for some $t_n>0$.
    We claim that $t_n\to\infty$.
    If it weren't so, indeed, there would be a bounded subsequence $t_{n_k}$ converting to some $\bar t$ and so, by continuity, $q_{n_k}=F^{t_{n_k}}(p_{t_k})\to F^{\bar t}(p)$.
    Hence $q$ would lie in the forward orbit of $p$, against the initial hypothesis.
    So $t_n\to\infty$ and therefore for each $\eps>0$ and $T>0$ we can find a $n$ such that $d(q,q_n)<\eps$, $d(p,p_n)<\eps$ and $t_n>T$, namely $q\in\Lambda^1_F(p)$.
\end{proof}
\begin{example}
    \label{ex: circle L^1}
    Continuing Example~\ref{ex: circle}, denote by $x$ any point in the interior of $\overset{\frown}{pq}$.
    Let us show that $\Lambda^1_F(p)=\overset{\frown}{pq}$.
    Indeed, $p\in\Lambda^1_F(p)$ because it is a fixed point and so its very orbit is an $(\eps,T_\eps)$-link from $p$ to itself for every $T_\eps$; $x\in\Lambda^1_F(p)$
    because, 
    by choosing $\eps>0$ small enough, it can take an arbitrarily large time to reach $x$ from a point in $\overset{\frown}{pq}$ within $\eps$ from $p$;
    finally, $q\in\Lambda^1_F(p)$ because from $x$ it takes an arbitrarily large time to get arbitrarily close to $q$.
    On the other side, $\Lambda^1_F(x)=\{q\}$ since it takes a finite time to get from $x$ to any other point $y$ in the interior of $\overset{\frown}{xq}$.
    Similarly happens in case of 
    $\overset{\frown}{qr}$ and $\overset{\frown}{rp}$.
    On the other side,
    $\Lambda^1_G(x)=\Bbb S^1$ for each $x\in\Bbb S^1$.
    Hence, while $\Lambda^1_G(\Lambda^1_G(x))=\Lambda^1_G(x)$ for each $x\in\Bbb S^1$, we have that $\Lambda^1_F(\Lambda^1_F(p))=\overset{\frown}{pr}$ and $\Lambda^1_F(\Lambda^1_F(x))=\overset{\frown}{qr}$ for each $x$ in the interior of $\overset{\frown}{qr}$.
    Notice that $\Lambda^1_F(\Lambda^1_F(\Lambda^1_F(p)))=\Bbb S^1$ and, more generally, $\Lambda^1_F(\Lambda^1_F(\Lambda^1_F(\Lambda^1_F(x)))=\Bbb S^1$ for each $x\in\Bbb S^1$.
\end{example}
As the example above shows, in general
$    \Lambda^1_F(\Lambda^1_F(x))\neq\Lambda^1_F(x).
$ 
Throughout the article we use the following notation:
$$
\Lambda^{1,k}_F(x)=\overbrace{\Lambda^1_F(\dots(\Lambda^1_F}^{k}(x))\dots),\;k=2,3,...\;.$$


\begin{example}
    \label{ex: [0,1]}
    Consider the flow on $X=[0,1]$ that keeps fixed all points 
    $$x_n=1-\frac{1}{n}, n\in\bN,$$ and moves rightwards every other point.
    Then 
    \begin{align*}
        &\Lambda^1_F\left(x_n\right)=\left[x_n,x_{n+1}\right], n\in\bN\\
        \text{and }&\\
        &\Lambda^1_F(x)=\left\{x_{n+1}\right\}\text{ for each }x\in\left(x_{n},x_{n+1}\right), n\in\bN.
    \end{align*}
    In particular, 
    $
    \Lambda^{1,k}_F(0)=\left[0,x_{k+1}\right]
    .$
\end{example}
Notice that, in the example above, the point 1 is not in $\Lambda^{1,k}(0)$ for any $k\in\bN$ but the right endpoints of $\Lambda^{1,k}(0)$ converge to 1 as $k\to\infty$.
This justifies the following definition.
\begin{definition}
    \label{def: L^2}
    The {\bf second prolongational limit set} of $p$ under $F$ is the set {\BF $\Lambda^2_F(p)$} of all points $q$ for which there are sequences $p_n,q_n\in X$ and $k_n\in\bN$ such that:
    \begin{enumerate}
        \item $p_n\to p$;
        \item $q_n\to q$;
        \item $q_n\in\Lambda^{1,k_n}_F(p_n)$.
    \end{enumerate}
    We define inductively the prolongational limit set of any given ordinal number $\alpha$ in the following way.
    Suppose that {$\Lambda^\beta_F(p)$} is defined for all $\beta<\alpha$ and let us use the notations $\Lambda^\beta_F(A)=\cup_{x\in A}\Lambda^\beta_F(x)$ and
    $$
    \Lambda^{\beta,k}_F(p)=\overbrace{\Lambda^\beta_F(\dots(\Lambda^\beta_F}^{k}(p))\dots),\;k=1,2,...\;.
    $$
    Then, if $\alpha$ is a successor ordinal, the {\BF $\alpha$-th prolongational limit set} of $p$ under $F$ is the set {\BF $\Lambda^\alpha_F(p)$} 
    of all points $q$ for which there are sequences $p_n,q_n\in X$, $k_n\in\bN$ such that:
    \begin{enumerate}
        \item $p_n\to p$;
        \item $q_n\to q$;
        \item $q_n\in\Lambda^{\alpha-1,k_n}_F(p_n)$.
    \end{enumerate}
    If $\alpha$ is a limit ordinal, the {\BF $\alpha$-th prolongational limit set} of $p$ under $F$ is the set {\BF $\Lambda^\alpha_F(p)$} 
    of all points $q$ for which there are sequences $p_n,q_n\in X$, $k_n\in\bN$ and ordinals $\beta_n<\alpha$ such that:
    \begin{enumerate}
        \item $p_n\to p$;
        \item $q_n\to q$;
        \item $q_n\in\Lambda^{\beta_n,k_n}_F(p_n)$.
    \end{enumerate}
    The {\BF second and $\alpha$-th prolongations} of $p$ under $F$ are defined analogously by replacing $\Lambda$ with $D$.
\end{definition}
\noindent
From the definition above it follows immediately that
$
\Lambda^{\alpha+1}_F(p) \supset \displaystyle\bigcup_{k\in\bN}\Lambda^{\alpha,k}_F(p).
$
\begin{proposition}[Auslander,1963~\cite{Aus63}]
    The following holds for every ordinal $\alpha$ and $p\in X$:
    \begin{enumerate}
        \item $D_F^\alpha(p) = \cO_F^+(p)\cup\Lambda^\alpha(p)$.
        \item $\Lambda_F^\alpha(p)$ (and so $D^\alpha_F(p)$) is invariant under $F$.
    \end{enumerate}
    
\end{proposition}
\begin{example}
    \label{ex: [0,1] bis}
    Continuing Example~\ref{ex: [0,1]}, now we see that, while $1\not\in\Lambda^{1,k}_F(0)$ for any $k\in\bN$, we have that $1\in\Lambda^2_F(0)$.
    Hence, $\Lambda^2_F(0)=[0,1]$ and
    $\Lambda^\alpha_F(0)=\Lambda^2_F(0)$ 
    for every ordinal $\alpha\geq2$.
    Similarly, for $x\in(x_{n},x_{n+1})$, $\Lambda^2_F(x)=[x_{n+1},1]$ and $\Lambda^\alpha_F(x)=\Lambda^2_F(x)$
    for every ordinal $\alpha\geq2$.
\end{example}
In the example above, we say that the $\Lambda^\alpha_F$ stabilize at 2.
The example below shows how to construct flows whose $\Lambda^\alpha_F$ stabilize at any countable ordinal.
\begin{example}
    Consider the flow $F$ in Examples~\ref{ex: [0,1]},\ref{ex: [0,1] bis} and now add fixed points so that, in each interval $[x_{n},x_{n+1}]$, we repeat the structure of the fixed points of $F$ in $[0,1]$.
    This means that in $[0,\frac{1}{2}]$ we add the fixed points $x_{1,n}=\frac{x_n}{2}$, in $[\frac{1}{2},\frac{2}{3}]$ the fixed points $x_{2,n}=\frac{1}{2}+\frac{x_n}{6}$ and so on.
    Every other point moves rightwards. 
    Denote this flow by $G$.
    By construction, loosely speaking, what was for $F$ a $\Lambda^1$ prolongational limit is for $G$ a $\Lambda^2$ prolongational limit and so on. 
    For instance:
    \begin{align*} 
    &\Lambda^{1,k}_G(0)=\left[0,\frac{x_{k+1}}{2}\right], k=1,2,\dots,\\
    \noalign{\medskip}
    &\Lambda^{2}_G(0)=\left[0,\frac{1}{2}\right],\\
    \noalign{\medskip}
    &\Lambda^{2,k}_G(0)=\left[0,x_{k+1}\right],\\
    \text{so that }&\\    &\Lambda^3_F(0)=[0,1].
    \end{align*}
    
    Since, for this flow, all prolongational limit sets are subsets of the ones at 0, it follow that $G$ stabilizes at 3, namely $\Lambda^\alpha_G(x)=\Lambda^3_G(x)$ for all $\alpha\geq3$ and $x\in[0,1]$.
    
    By repeating this trick any finite number of times, for each $k\geq 3$ one can get a flow $G_k$ that stabilizes at $k$ and such that 
    $\Lambda^k_F(0)=[0,1]$.
    Now one can use these flows as follows. 
    Set a flow $G_2$ on $[0,\frac{1}{2}]$, a flow $G_3$ on $[\frac{1}{2},\frac{2}{3}]$ and so on. 
    As usual, every other point moves rightwards.
    Denote this flow by $G_\omega$, where $\omega$ is the first limit ordinal.
    Hence, $\Lambda^2_{G_\omega}(0)=[0,\frac{1}{2}]$, $\Lambda^{3,2}_{G_\omega}(0)=[0,\frac{2}{3}]$ and, in general, $\Lambda^{n+1,n}_{G_\omega}(0)=[0,x_{n+1}]$.
    Then, 
    $\Lambda^\omega_{G_\omega}(0)=[0,1]$ so that, by the same arguments used in the previous cases, $G_\omega$ stabilizes at $\omega$.
    If now one extends $G_\omega$ to a flow on $[0,k+1]$ so that, on each interval $[i,i+1]$, the flow restricts to $G_2$, then the resulting flow stabilizes at $\omega+k$.
    If $G_\omega$ is extended to a flow on $[0,2]$ so that its restriction to $[1,2]$ is again $G_\omega$, then the resulting flow stabilizes at $2\omega$.
    If $G_\omega$ is extended to $[0,2]$ so that 2 is a fixed point and on each $[2-\frac{1}{i},2-\frac{1}{i+1}]$, $i=1,2,\dots$, the flow is $G_\omega$, then this flow stabilizes at $\omega^2$ and so on.
\end{example}
\begin{definition}
    \label{def: rank}
    We denote by {\BF $\Lambda_F^\star(p)$} the largest of the $\Lambda_F^\alpha(p)$.
    We say $F$ is a {\BF flow of rank $\alpha$} if $\alpha$ is the smallest ordinal such that $\Lambda_F^\star(p)=\Lambda_F^\alpha(p)$ for every $p\in X$. 
    Analogously, we denote by {\BF $D^\star_F(p)$} the set $\cO_F^+(p)\cup\Lambda_F^\star(p)$.
\end{definition}
It turns out that are no flows of rank higher than $\omega_1$, the first uncountable ordinal:
%
\begin{thmX}[Auslander and Seibert, 1964~\cite{AS64}]
    \label{thm: Lambdas}
    For any flow $F$,
    $\Lambda^{\omega_1}_F=\Lambda^\star_F$. 
\end{thmX}

We will show in Theorem~\ref{thm1} that every regular flow in the plane has rank either 0, 1 or 2 and that $\Lambda_F^2(p)=\cup_{n\in\bN}\Lambda^{1,n}_F(p)$ for every $p\in X$.
%
\begin{remark}
    The structure of the $\Lambda^\alpha$ highlighted by Theorem~\ref{thm: Lambdas} is somehow reminiscent of the Cantor-Bendixon rank of a closed set.
    In that case, one starts with a closed set $C$ and sets $C^1=C$, $C^\alpha=(C^{\alpha-1})'$ (the Cantor-Bendixon derivative of the set, namely the set minus its isolated points) if $\alpha$ is a successor ordinal and $C^\alpha=\cap_{\beta<\alpha}C^{\beta}$ otherwise.
    Then, it can be proved that, for every $C$, there exists an ordinal $\gamma$ (called rank of $C$), smaller than the first uncountable ordinal,  such that $C^{\gamma+1}=C^\gamma$.
    Then, for every other $\alpha>\gamma$, $C^\alpha=C^{\gamma}$.
\end{remark}
\begin{definition} 
    We say that $p$ is a {\bf non-wandering point} for $F$ if $p\in\Lambda^1_F(p)$.
    We denote by {\BF $NW_F$} the set of all non-wandering points of $F$.
\end{definition}
\begin{proposition}
    \label{prop: oa}
    For every $p\in X$, $\omega_F(p)\cup\alpha_F(p)\subset NW_F$.
\end{proposition}
\begin{proof}
    Let $q\in\omega_F(p)$.
    Then there is a sequence $t_n\to\infty$ such that $F^{t_n}(p)\to q$.
    Hence, for every $\eps>0$ there is a $k>0$ such that $d(q,F^{t_k}(p))<\eps$ and, for every $T>0$, there is a $k'>k$ such that $t_{k'}\geq T$ and $d(q,F^{t_{k'}}(p))<\eps$.
    In other words, for every $\eps>0$ and $T>0$ there is an $(\eps,T)$-link from $q$ to itself, namely $q\in\Lambda^1_F(q)$, i.e. $q\in NW$.
    The same argument shows that $\alpha_F(p)\subset NW_F$.
\end{proof}
\begin{proposition}
    \label{prop: OL}
    For every $p\in X$, $\cO^+_F(p)\cap\Lambda^1_F(p)\subset NW_F$.
\end{proposition}
\begin{proof}
    Let $q\in\cO^+_F(p)\cap\Lambda^1_F(p)$.
    Then:
    \begin{enumerate}
        \item since $q\in\cO^+_F(p)$, there is $t_0\geq0$ such that $q=F^{t_0}(p)$;
        \item since $q\in\Lambda^1_F(p)$, for every $\eta>0$ there is a $z$ within $\eta$ from $p$ and a $t>0$ such that $d(q,F^t(z))<\eta$ and $t\to\infty$ as $\eta\to0$;
        \item by the continuity of $F^{t_0}$, for every $\eps>0$ there is $\eta>0$ such that $d(z,p)<\eta$ implies $d(F^{t_0}z,q)<\eps$.
    \end{enumerate}
    We claim that $q\in NW_F$.
    Indeed, let $\eps>0$, choose $\eta>0$ as in (3) above (if $\eta>\eps$, set $\eta=\eps$) and take a corresponding $z$ as in (2) above.
    The orbit through $z$ passes within $\eps$ from $q$, so there exists a $z'$ within $\eps$ from $q$ and a $t'>0$ such that $d(q,F^{t'}(z'))<\eps$. 
    Hence, $q\in NW_F$.
\end{proof}
\begin{proposition}
    \label{prop: NW_G}
    Let $G$ be the inverse flow of $F$, namely $G^t(p) = F^{-t}(p)$.
    Then:
    \begin{enumerate}
        \item $p\NWto q$ if and only if $q\NWGto p$;
        \item $NW_F=NW_G$.
    \end{enumerate}
\end{proposition}
\begin{proposition}
    \label{prop:NW}
    For every flow $F$, $NW_F$ is closed and $F$-invariant.
\end{proposition}
\begin{proof}
    Let $x_n$ be a sequence of non-wandering points of $F$ converging to some $\bar x\in X$.
    Note that, given any point $y$ closer to $\bar x$ than $\eps$, every $(\eps,T)$-link from $y$ to itself is also a $(2\eps,T)$-link from $\bar x$ to itself. 
    Hence, since there are points of $NW_F$ arbitrarily close to $\bar x$, this implies that $NW_F$ is closed.
    Now, let $x\in NW_F$. 
    We claim that, as a consequence of the continuity of $F$, $\cO^+_F(x)\subset NW_F$.
    Indeed, given any $\tau>0$, for every $\eps>0$ there is $\eta>0$ (depending on $\tau$) such that $d(x,y)<\eta$ implies that $d(F^\tau x,F^\tau y)<\eps$.
    We can assume without loss of generality that $\eta\leq\eps$.
    Hence, given a $(\eta,T)$-link $\ell$ from $x$ to itself starting at $y$ and ending at $F^Ty$, so that $d(x,y)<\eta$ and $d(x,F^Ty)<\eta$, we have that $F^\tau\circ\ell$ is a $(\eps,T)$-link from $F^\tau x$ to itself starting at $F^\tau y$ and ending at $F^{T+\tau} y$.
    Since such a link exists for every $\eps>0,T>0$, it follows that $F^\tau x\in NW_F$.
    Finally, since $F^\tau$ is invertable for every $\tau$, $NW_F$ also contains the backward orbit of each point and so it is invariant.
\end{proof}

\noindent
{\bf Recurrence and Lyapunov functions.} It turns out that there is a strong relation between non-wandering points and Lyapunov functions, as defined below.
\begin{definition}
    \label{def: Lyap}
    We say that $L\in C^0(X)$ is a {\bf Lyapunov function} for $F$ if, for all $x\in X$, $y\in\cO_F^+(x)$ implies $L(x)\geq L(y)$, namely $L$ is non-increasing on each orbit.
    We say that $L$ is a {\bf strict Lyapunov function} for $F$ if $y\in\cO_F^+(x)$, $y\neq x$, implies $L(x)> L(y)$. 
    We denote by $\cL_F(X)$ the set of Lyapunov functions of $F$.
\end{definition}
%
%
%
\begin{definition}
    \label{def: stream}
    We say that a binary relation $D$ on $X$ is a {\bf quasi-order} if it is reflexive and transitive.
    We say that a quasi-order $E$ is an {\bf extension} of $D$ if $E\supset D$.
    We say that $D$ is a {\BF $F$-stream} (or simply a {\bf stream}, when there is no ambiguity) if $D$ is a transitive and topologically closed extension of $\NW_F$.
\end{definition}
\begin{proposition}
    \label{prop: A_F}
    The relation
    $$\cA_F=\{(x,y):L(x)\geq L(y) \text{ for all } L\in\cL_F(X)\}$$
    is the smallest $F$-stream.
\end{proposition}
\begin{proof}
    Let us prove first that $\cA_F$ is a $F$-stream.
    First of all, $\cA_F$ is closed, since it is defined by a closed relation on continuous functions.
    Moreover, $\cO_F\subset\cA_F$ by the very definition of $\cA_F$.
    Hence, $\cA_F\supset\overline{\cO_F}=\NW_F$.
    The reader can verify that $\cA_F$ is both reflexive and transitive, so it is a closed extension of $\NW_F$.
    
    Finally, suppose that there is a stream $D$ strictly smaller than $\cA_F$.
    Then there must be at least a pair $(x,y)$ such that $L(x)\geq L(y)$ for all $L\in\cL_F(X)$ and $(x,y)\not\in D$.
    By Lemma~A in~\cite{DLY24}, there is a function $K\in C^0(X)$ such that $K(x')\geq K(y')$ for each $(x',y')\in D$ and $K(x)<K(y)$.
    This leads to a contradiction: since $D\supset\cO_F$, then $K\in\cL_F(X)$, but then, by the very definition of $\cA_F$, $(x,y)$ cannot belong to $\cA_F$, against the initial hypothesis.
\end{proof}
The smallest stream was first considered by the first author in~\cite{Aus63}, which is why we also call it the Auslander stream.
Notice that, since $\cA_F$ is a transitive extension of $\NW_F$ then, 
when $\NW_F$ is transitive, $\cA_F=\NW_F$.
\begin{example}
    Let $F$ be the flow of the ODE $x'=|\sin(\pi x)|$ on $X=\bR$.
    Denote by $n_x$ the smallest integer with $n_x>x$.
    Then, under $F$, every point with integer coordinate is fixed while every other point $x$ moves monotonically rightwards and asymptotes to $n_x$.
    Then $\cO_F(x)=[x,n_x)$ and, for $\eps>0$ small enough, any $(\eps,T)$-link from $x$ to $y$ satisfies the following relation:
    $n_x-1<x-\eps<F^T(x)+\eps<n_x+\eps.$
    Hence, $\Down_{\NW_F}(x)=[x,n_x]$.
    Since $\cA_F$ is transitive and contains $\NW_F$, it follows immediately that $\Down_{\cA_F}(x)=[x,\infty)$.
\end{example}
\noindent
{\BF Prolongational sets and $\cA_F$.}
Next result shows the relation between the smallest stream and the prolongational relation.
\begin{thmX}[Auslander, 1963~\cite{Aus63}]
    \label{thm: Down A}
    $\Down_{\cA_F}(p)=D_F^\star(p)$.
\end{thmX}
\begin{definition}
    \label{def: gener rec pts}
    We write $x\Ato y$ if $(x,y)\in\cA_F$.
    We say that $x$ is a {\bf generalized recurrent point} if either $x$ is fixed or there is a $y\neq x$ such that $x\Ato y$ and $y\Ato x$.
    We denote by $A_F$ the set of all generalized recurrent points of $F$.
\end{definition}
Notice that, as $\NW_F\subset\cA_F$, we have also that $NW_F\subset A_F$.
\begin{proposition}[Auslander, 1963~\cite{Aus63}]
    The following are equivalent:
    \begin{enumerate}
        \item $x\in A_F$.
        \item $x\in\Lambda^\alpha_F(x)$ for some ordinal $\alpha$. 
        \item $\cO^-_F(x)\subset D^\alpha_F(x)$ for some ordinal $\alpha$.
        \item $F^t(x)\in D^\alpha_F(x)$ for some ordinal $\alpha$.
        \item $D^\alpha_F(x)=\Lambda_F^\alpha(x)$ for some ordinal $\alpha$.
    \end{enumerate}
\end{proposition}
The following result gives a characterizations of $A_F$ in terms of Lyapunov functions.
%
\begin{thmX}[Auslander, 1963 (Thm. 2 in ~\cite{Aus63})]
    \label{prop: lyap2}
    There is a function $L\in \cL_F(X)$ such that:
    \begin{enumerate}
        \item $p\in A_F$ if and only if $L$ is constant on $\cO_F(p)$;
        \item if $p\not\in A_F$, then $L(p)>L(F^t(p))$ for each $t>0$.
    \end{enumerate}
\end{thmX}
%

%
\black\medskip\noindent
{\bf Wandering flows.}
The focus of this article is on flows for which every point is wandering, as defined below.
%
\begin{definition}
    We say that $F$ is {\bf wandering} if $NW_F=\emptyset$.
\end{definition}
Recall that every recurrent point is non-wandering. 
Hence, given a wandering flow $F$, for any sequence $x_n\to x$ and $t_n\to\infty$, there cannot be any converging subsequence of $F^{t_n}(x_n)$ (if any) that converges to $x$.
When $\omega(x)\neq\emptyset$, for every sequence $t_n\to\infty$ the sequence $F^{t_n}(x)$ has converging subsequences (to some point of $\omega(x)$).
When $F$ is wandering, though, $\omega(x)=\emptyset$ for each $x$ and 
we expect that, given generic sequences $t_n\to\infty$ and $x_n\to x$, the sequence $F^{t_n}(x_n)$ has no converging subsequence. 
For instance the example below shows that, even when $F^{t_n}(x_n)$ converges (to some point in $\Lambda^1_F(x)$), there are sequences $\tau_n$ with $\tau_n\to\infty$ and $\tau_n< t_n$ such that no subsequence of $F^{\tau_n}(x_n)$ does converge.
%
%
\begin{example}
    Consider the complete vector field $\xi(x,y)=(2y,1-y^2)$ and set $p=(0,-1)$, $q=(0,1)$.
    Denote by $F_\xi$ the flow of $\xi$. Then $q\in\Lambda^1_{F_\xi}(p)$.
    Indeed, let $p_n=(0,\frac{1}{n}-1)$ and $t_n=\ln(2n-1)$.
    Then 
    $$
    F_\xi^t(p_n) 
    = 
    \left(2\ln\frac{e^{2t}+2n-1}{2n}-2t,\frac{e^{2t}-(2n-1)}{e^{2t}+2n-1}\right)
    $$
    and, for $n\to\infty$, we have that $p_n\to p$, $t_n\to\infty$ and 
    $F_\xi^{t_n}(p_n)=(0,1-\frac{1}{n})\to q$.
    Now, set $\tau_n=\frac{1}{2}\ln(2n-1)$.
    Then $\tau_n\to\infty$ and $\tau_n<t_n$ but the sequence
    $$
    F_\xi^{\tau_n}(p_n) 
    = 
    \left(\ln\frac{2n-1}{n^2},0\right)
    $$
    has no converging subsequence.
\end{example}
\black
\black

We show below that being wandering is inherited by inverses, factors and products.
\begin{definition}
    A flow $G$ on a topological space $Y$ is a {\bf factor} of $F$ if there is a continuous, surjective and equivariant map $\phi:X\to Y$, namely $\phi(F^t(x))=G^t(\phi(x))$ for every $t\in\bR$.
    The flow $F$ is called an {\bf extension} of $G$.
\end{definition}
\begin{proposition}
    The following holds:
    \begin{enumerate}
        \item The inverse of 
        a wandering flow 
        is wandering.
        \item Any extension of 
        a wandering flow 
        is wandering.
        \item The product of wandering flows is wandering.
    \end{enumerate}
\end{proposition}
\begin{proof}
    (1) The non-wandering set $NW_F$ is invariant under $F$ and therefore also under $F^{-1}$, so $NW_F=NW_{F^{-1}}$.

    (2) Assume that the flow $G$ on $Y$ is a factor of $F$. 
    Then, if $x\in NW_F$, we have that $\phi(x)\in NW_G$.
    Hence, if $G$ is wandering, necessarily $F$ must be wandering.

    (3) Both $F$ and $G$ are factors of $F\times G$, so this is a corollary of (2).
\end{proof}
A factor of a wandering flow is not necessarily wandering.
For instance, any factor $G$ acting on a compact $Y$ has a non-empty non-wandering set, whether its extension $F$ is wandering or not.
Below we present a necessary condition for $G$ to be wandering.
\begin{proposition}
    Let the flow $G$ on $Y$ be a factor of $F$ with equivariant map $\phi:X\to Y$.
    A necessary condition for $G$ to be wandering is that, if $x\in X$ and $t_n\to\infty$ is such that $F^{t_n}(x)\to x'$, then $\phi(x)\neq\phi(x')$.
\end{proposition}
The following relations between orbits, limit sets and prolongational limit sets are fundamental for our article. 
\begin{proposition}
    Let $F$ be a wandering flow. 
    Then, for all $p$:
    \begin{enumerate}
        \item $\cO^+_F(p)\cap\Lambda^1_F(p)=\emptyset$;
        \item $\omega_F(p)=\alpha_F(p)=\emptyset$.
    \end{enumerate}
\end{proposition}
\begin{proof}
    Case (1) is an immediate consequence of Proposition~\ref{prop: Down}.
    Case (2) is an immediate consequence of Propostion~\ref{prop: oa}.
\end{proof}
\begin{corollary}
    \label{cor: open orbits}
    Each orbit of a wandering flow is unbounded for both $t\to+\infty$ and $t\to-\infty$.
\end{corollary}
%
%
\section{Regular planar flows}
\label{sec: wandering}
Our ultimate goal is classifying wandering flows according to the following equivalence relation.
\begin{definition}
    Two flows $F,G$ on $X$ are {\bf locally topologically equivalent} if, for every $p$, there is an open neighborhood $U_p\subset X$ of $p$ and a homeomorphism $\phi_p:U_p\to X$ such that, for every $q\in U_p$, 
    $\phi_p(\cO_F(q)\cap U_p)\subset\cO_G(\phi_p(q))$. 
    We say that $F$ and $G$ are {\bf topologically equivalent} if
    there exists a homeomorphism of $X$ that sends orbits of $F$ in orbits of $G$ and preserves the orientation of the orbits.
\end{definition}
\noindent
{\bf Regular flows.} A full topological classification of flows is in general highly non-trivial unless one poses some suitable restriction.
In this article, we consider the case of {\em regular} flows as defined below.
\begin{definition}
    We say that $F$ is {\bf regular} if it is locally topologically equivalent to the flow of a constant vector field. 
\end{definition}
\begin{example}
    The flow $F$ of a $C^1$ vector field $\xi$ without zeros is $C^1$ regular.
    This is an immediate consequence of the Rectification Theorem (e.g. see~\cite{Arn92b}), stating that, in the neighborhood of any point where $\xi(p)\neq0$, there is a $C^1$ diffeomorphism mapping $\xi$ into a constant vector field.
    Given any homeomorphism of the plane $\phi$, the flow $F_\phi^t=\phi^{-1}\circ F^t\circ\phi$ is $C^0$ regular.
\end{example}
\begin{proposition}
    Let $F$ be a regular planar flow.
    Then $F$ is wandering.
\end{proposition}
\begin{proof}
    By hypothesis, $F$ cannot have fixed points or periodic orbits.
    By the Poincar\'e-Bendixon theorem, this means that $NW_F=\emptyset$. 
\end{proof}
\medskip\noindent
{\bf Regular flows and foliations.}
By Corollary~\ref{cor: open orbits}, each orbit of a regular flow is unbounded for $t\to\pm\infty$.
Hence, each regular flow determines a family of continuous open curves that go to infinity at both ends and whose union is the whole plane.
\begin{definition}
    \label{def: fol}
    By {\bf foliation of the plane} we mean a family of continuous planar curves with the following properties:
    \begin{enumerate}
        \item for every $\alpha\in\cF$, $\alpha$ goes to infinity at both ends; 
        \item for every $\alpha,\beta\in\cF$, either $\alpha=\beta$ of $\alpha\cap\beta=\emptyset$;
        \item $\displaystyle\bigcup_{\alpha\in\cF}\alpha=\bR^2$;
        \item $\cF$ is locally homeomorphic to a family of parallel lines.
    \end{enumerate}
\end{definition}
\begin{proposition}
    The set of orbits of any regular flow $F$ on the plane defines a foliation of the plane.
\end{proposition}
\begin{proof}
    We need to show that the set of all orbits of $F$ satisfies the three points of Definition~\ref{def: fol}.
    Points (1) and (2) hold by the definition of flow.
    Point (3) follows from the regularity of $F$.
\end{proof}
%
%
%
%
Next example shows the importance of requiring flows to be regular.
%
\begin{figure}
    \centering
    \includegraphics[width=0.6\linewidth]{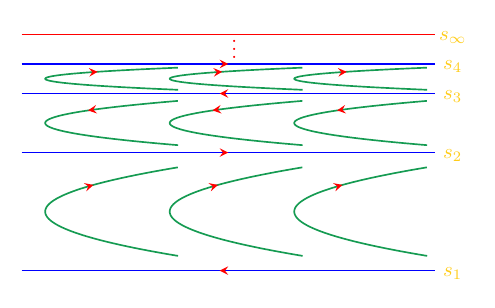}
    \caption{An example of a family of continuous mutually disjoint curves filling the plane. 
    No point on $s_\infty$ has a neighborhood where the family is topologically equivalent to a family of parallel straight lines.
    }
    \label{fig: nr}
\end{figure}
\begin{example}
   Consider the family of curves $\cF$ in Figure~\ref{fig: nr}.
   Outside of $s_\infty$, a consistent orientation can be given to each curve of $\cF$ so that each of them is the orbit of a flow.
   On the other side, any vertical line segment through any point of $s_\infty$ meets all curves $s_k$ for $k$ large enough.
   Hence, arbitrarily close to $s_\infty$, there are orbits where points move in opposite directions.
   This implies that $s_\infty$ cannot be the orbit of a flow: every semiflow $F$ having the curves in $\cF\setminus\{s_\infty\}$ as orbits, necessarily has each point of $s_\infty$ as a fixed point.
   In particular, $\cF$ is not the family of orbits of a regular flow.
\end{example}

\medskip\noindent
{\bf Inseparable orbits and separatrices.}
The quotient topology of the space of orbits plays a fundamental role in the classification of flows. 
This topology is non-Hausdorff except for the trivial case of flows topologically equivalent to the flow of a constant vector field,
which justifies the following definition.
%
\begin{definition}
    \label{def: inseparable}
    A set $A\subset X$ is {\BF saturated} if $p\in A$ implies that $\cO_F(p)\subset A$.
    We say that an orbit $s_1$ of $F$ is {\bf inseparable} from another orbit of $F$ $s_2\neq s_1$ if, for every saturated neighborhoods $U_1$ of $s_1$ and $U_2$ of $s_2$ we have that $U_1\cap U_2\neq\emptyset$.
    In this case, we say that $s_1$ and $s_2$ are {\bf separatrices} for the flow $F$.
    We denote by $\cS_F$ the set of all separatrices of $F$.
\end{definition}
As we show below, the notion of inseparability is strictly related to the Auslander stream.
\begin{proposition}
    A separatrix $s_1\in\cS_F$ is inseparable from another separatrix $s_2\in\cS_F$ if and only if, possibly after relabeling the two separatrices, for every $p\in s_1$ and every $q\in s_2$ we have that $p\NWto q$.
\end{proposition}
\begin{proof}
    Assume that $s_1$ is inseparable from $s_2\neq s_1$ and let $p\in s_1$ and every $q\in s_2$.
    Let $J_p$ (resp. $J_q$) be a segment passing through $p$ and such that every orbit of $F$ meets $J_p$ (resp. $J_q$) in at most one point -- such segment always exists because $F$ is regular -- and denote by $U$ (resp. $V$) the set of all orbits that intersect $U$ (resp. $V$).
    Since $s_1$ and $s_2$ are inseparable, $s_1\cap s_2\neq\emptyset$.
    Hence, for every $\eps>0$ there is an orbit of $F$ that passes within $\eps$ from $s_1$ and within $\eps$ from $s_2$.
    Hence, $p\NWto q$.
\end{proof}
Notice that, for points $p,q$ belonging to distinct orbits, $p\NWto q$ is equivalent to $q\in\Lambda_F^1(p)$
%
    and that the relation of ``being inseparable from'' is symmetric but not transitive (e.g. see Ex.~\ref{ex:four seps}).  
%
\begin{figure}
    \centering
    \includegraphics[width=0.6\linewidth]{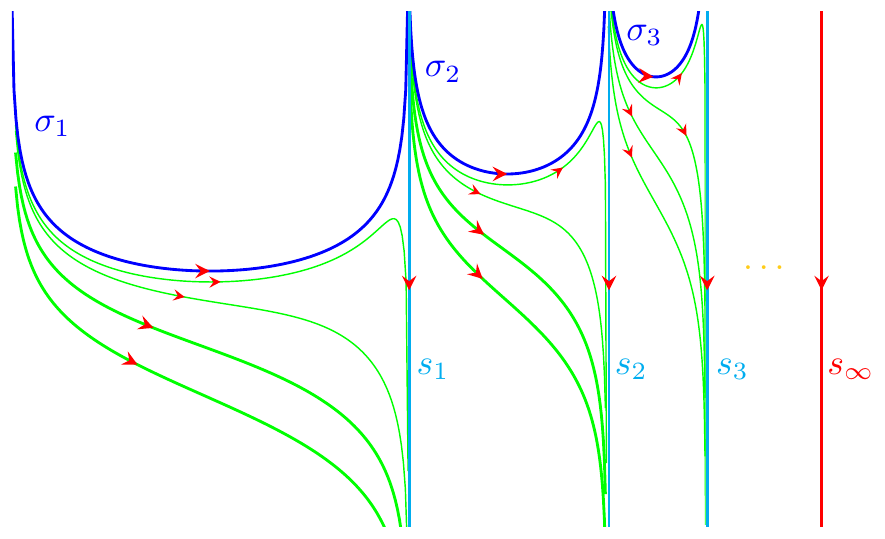}
    \caption{An example of a flow with a non-closed set of separatrices. 
    Separatrices are painted in blue, cyan and red.
    Orbits that are not separatrices are painted in green.
    }
    \label{fig: waz}
\end{figure}

As shown in the example below, the set $\cS_F$ is not necessarily closed. 
\begin{example}
    Consider the flow $F$ whose orbits are sketched in Figure~\ref{fig: waz}.
    Then $\cS_F=\cup_{k=1}^\infty \sigma_k\cup_{k=1}^\infty s_k\cup s_\infty$.
    Each pair $s_k,\sigma_k$ is a pair of inseparable separatrices. 
    The orbit $s_\infty$ is not a separatrix but it the limit set of the separatrices $s_k$ as $k\to\infty$.
    The orbits $\sigma_k$ are diverging to infinity as $k\to\infty$.
    It is possible to exploit the idea illustrated in Figure~\ref{fig: waz} to produce examples of flows whose separatrices accumulate on any finite or infinite countable number of orbits, including the case of orbits that are dense in the plane.
    The first example of a planar flow with separatrices dense on the plane was found 
    by Wazewski~\cite{Waz34}
    (see also~\cite{HR57} and~\cite{DeL14c}).
\end{example}
\begin{definition}
    \label{def: inseparable}
    Let $s_1,s_2$ be distinct inseparable separatrices and let $p\in s_1$ and $q\in s_2$.
    If $p\NWto q$
    we say that {\BF $s_1$ precedes $s_2$} and we write {\BF $s_1\prec s_2$}, and viceversa if $q\NWto p$.
\end{definition}
The proposition below grants that the definition above does not depend on the choice of $p$ and $q$.
\begin{proposition}
    \label{prop: s_1 prec s_2}
    $s_1\prec s_2$ if and only if 
    $z\NWto w$ for every $z\in s_1$ and $w\in s_2$.
\end{proposition}
\begin{proof}
    Assume that $s_1\prec s_2$.
    Then, by definition, there are $p\in s_1$ and $q\in s_2$ such that $p\NWto q$.
    Let $z\in s_1$ and $w\in s_2$.
    There are four cases depending on the relative positions of $z$ and $p$ and of $q$ and $w$.
    Given any $\eps>0$, there are neighborhoods $U_p,U_q,U_z,U_w$ such that every orbit in $U_\alpha$ is within $\eps$ from $\alpha$ for each $\alpha=p,q,z,w$.
    Then $U_p\cap U_z$ is a neighborhood of $s_1$ and $U_q\cap U_w$ is a neighborhood of $s_2$.
    Since $s_1$ and $s_2$ are inseparable, $U_p\cap U_z\cap U_q\cap U_w$ is non-empty, namely there is an orbit that passes within $\eps$ from both $z$ and $w$.
    Hence, $z\NWto w$.
\end{proof}
%
\black
\begin{corollary}
    If $\cO_F(p)$ is not a separatrix, then $\Lambda_F^1(p)=\emptyset$.
    If $\cO_F(p)$ is a separatrix, then $\Lambda_F^1(p)$ is the union of all points belonging to all separatrices $s$ of $F$ such that $\cO_F(p)\prec s$.
\end{corollary}
Notice that the first prolongation of a point on a separatrix $s$ could be empty as well, in case $s$ does not precede any separatrix.

\begin{proposition}
    \label{lemma: L^(1,k)}
    Let $p,q\in X$. Then $q\in\Lambda^{1,k}_F(p)$ if and only if there are separatrices $s_1,\dots,s_{k+1}$ of $F$ such that:
    \begin{enumerate}
        \item $p\in s_1$;
        \item $q\in s_{k+1}$;
        \item $s_{i}\prec s_{i+1}$ for all $i=1,\dots,k$.
    \end{enumerate}
\end{proposition}
\begin{proof}
    For $k=1$, this is just Def.~\ref{def: inseparable}.
    Let now $k=2$, namely assume $q\in\Lambda^{1,2}(p)$.
    Then there exists $r\in X$ such that $r\in\Lambda^1_F(p)$ and $q\in\Lambda^{1}_F(r)$.
    Hence, 
    $\cO_F(p),\cO_F(r),\cO_F(q)$ are separatrices and $\cO_F(p)\prec\cO_F(r)$, $\cO_F(r)\prec\cO_F(q)$.
    The claim follows by induction.
\end{proof}
In the following examples we evaluate $\Lambda^\star$ and rank (see Definition~\ref{def: rank}) of several smooth wandering flows.
Throughout the examples, will denote by $\fL_\xi H$ the directional (Lie) derivative of the function $H$ along the vector field $\xi$, namely 
$$
\fL_\xi H(p) = {d\over dt}\bigg|_{t=0}(H\circ F^t_\xi)(p),
$$
where $F^t_\xi$ is the (local) flow of $\xi$.
Notice that, if $K\in C^1(X)$ is a strict Lyapunov function for the flow of a continuous vector field $\xi$, then $\fL_\xi K>0$ at every point.

\begin{figure}
    \centering
    \includegraphics[width=5cm]{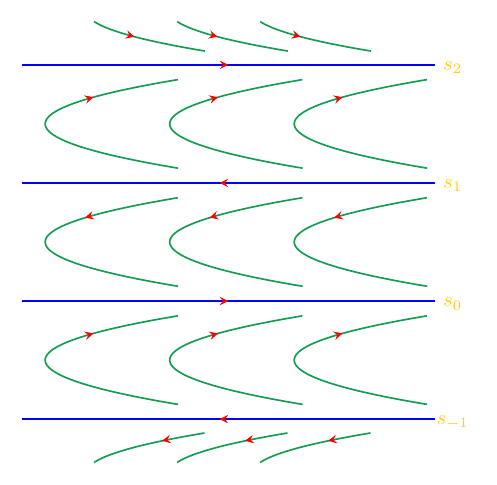}
    \hspace{0.5cm}
    \includegraphics[width=5cm]{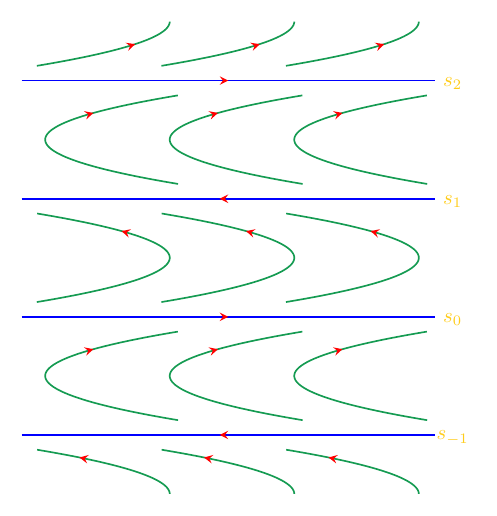}
    \caption{(left) Sketches of the flows of the vector fields $\xi(x,y)=(\sin y,\cos y)$ (left) and $\xi(x,y)=(\sin y,\cos^2 y)$ (right).
    Their separatrices are painted in blue.
    }
    \label{fig: infinitely many separatrices}
\end{figure}

\black
\begin{example}[\bf A trivial example]
    \label{ex:trivial}
    Consider the flow $F$ of the vector field $\xi(x,y)=(1,0)$ on $\bR^2$.
    Notice that $\fL_\xi H(x,y)=0$ for $H(x,y)=y$, corresponding to the fact that each orbit of $F$ is a horizontal straight line.
    The function $K(x,y)=x$ is a $C^\omega$ strict Lyapunov function for $F$.
    In this case, there are no separatrices, so that:
    \begin{enumerate}
        \item $\Lambda^\star_F(p)=\emptyset$ for every $p$;
        \item $rank(F)=0$;
        \item $\NW_F=\cA_F=\cO_F$;
        \item $NW_F=A_F=\emptyset$.
    \end{enumerate}
\end{example}
\begin{figure}
    \centering
    \includegraphics[width=6cm]{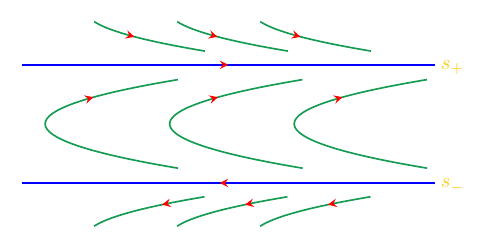}
    \hspace{0.5cm}
    \includegraphics[width=4.5cm]{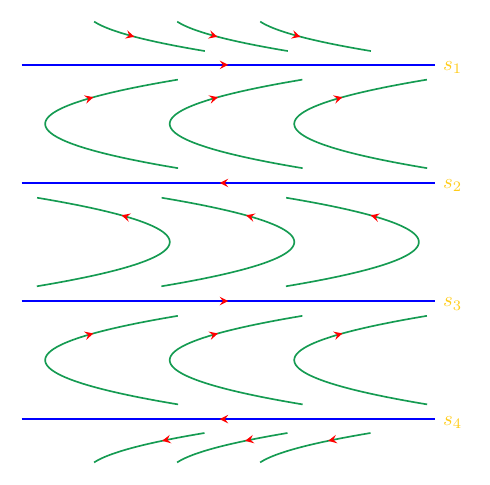}
    \caption{(left) Sketch of the flow $F$ of the vector field $\xi(x,y)=(2y,1-y^2)$.
    The two separatrices of $F$ are painted in blue.
    (right) Sketch of a flow $F$ with four separatrices, painted in blue.
    }
    \label{fig: finitely many seps}
\end{figure}
\begin{example}[\bf A case with two separatrices]
    \label{ex:two seps}
    Let $F$ be the $C^\omega$ flow of the $C^\omega$ vector field $\xi(x,y)=(2y,1-y^2)$ on $\bR^2$, sketched in Fig.~\ref{fig: finitely many seps}(left).
    
    Notice that $\fL_\xi H(x,y)=0$ for $H(x,y)=(y^2-1)e^x$, so that 
    $$
    \cO_F(x_0,y_0)\subset
    \{(x,y): H(x,y)=H(x_0,y_0)\}.
    $$
    This flow has precisely two separatrices, namely the lines
    $s_\pm = \{y=\pm1\},$ with $s_-\prec s_+$.
    Hence
    $$
    \Lambda^1_F(p)=
    \begin{cases}
        \emptyset,&p\not\in s_-;\\
        s_+,&p\in s_-,
    \end{cases}
    $$
    so that: 
    \begin{enumerate}
        \item $\Lambda^\star_F(p)=\Lambda^1_F(p)$ for each 
        $p\in\bR^2$;
        \item $rank(F)=1$;
        \item $\NW_F=\cA_F$;
        \item $NW_F=A_F=\emptyset$.
    \end{enumerate}
    %
    %
\end{example}
\begin{example}[\bf A case with four separatrices]
    \label{ex:four seps}
    Let $F$ be the planar flow whose orbits are sketched in Fig.~\ref{fig: finitely many seps}(right).
    In this case there are four separatrices $s_1,s_2,s_3,s_4$ and they are such that
    $s_1\prec s_2$, $s_2\prec s_3$ and $s_3\prec s_4$.
    Notice that $\prec$ is not transitive and those above are the only relations between the separatrices. 
    In particular, $s_1\not\prec s_3$ and so on.
    Hence
    $$
    \Lambda^1_F(p)=
    \begin{cases}
        \emptyset,&p\not\in s_1\cup s_2\cup s_3;\\
        s_2,&p\in s_1;\\
        s_3,&p\in s_2;\\
        s_4,&p\in s_3
    \end{cases}
    $$
    and
    $$
    \Lambda^2_F(p)=
    \begin{cases}
        \emptyset,&p\not\in s_1\cup s_2\cup s_3;\\
        s_2\cup s_3\cup s_4,&p\in s_1;\\
        s_2\cup s_3,&p\in s_2;\\
        s_2,&p\in s_3,
    \end{cases}
    $$
    so that: 
    \begin{enumerate}
        \item $\Lambda^\star_F(p)=\Lambda^2_F(p)$ for each
        $p\in\bR^2$;
        \item $rank(F)=2$;
        \item $\NW_F\subsetneq\cA_F$;
        \item $NW_F=A_F=\emptyset$.
    \end{enumerate}
\end{example}
\begin{example}[\bf A first case with infinitely many separatrices]
    \label{ex:inf seps 1}
    Let $F$ be the flow of the $C^\omega$ vector field $\xi(x,y)=(\sin y,\cos y)$, whose orbits are sketched in Fig.~\ref{fig: infinitely many separatrices}(left).
    
    Notice that $\fL_\xi H(x,y)=0$ for $H(x,y)=e^x\cos y$, so that 
    $$
    \cO_F(x_0,y_0)\subset
    \{(x,y): H(x,y)=H(x_0,y_0)\}.
    $$
    The separatrices of this flow are the horizontal lines $s_k=\{y=k\pi+\pi/2\}$, $k\in\bZ$.
    The relations between the separatrices are the following:
    $$
    s_{2n+1}\prec s_{2n},\;
    s_{2n+1}\prec s_{2(n+1)}
    $$
    for every $n\in\bZ$.
    Hence 
    $$
    \Lambda^1_F(p)=
    \begin{cases}
        \emptyset,&p\not\in s_{2n+1}\text{ for any }n\in\bZ;\\
        s_{2(n+1)}\cup s_{2n},&p\in s_{2n+1}\text{ for some }n\in\bZ,\\
    \end{cases}
    $$
    so that: 
    \begin{enumerate}
        \item $\Lambda^\star_F(p)=\Lambda^1_F(p)$ for each 
        $p\in\bR^2$;
        \item $rank(F)=1$;
        \item $\NW_F=\cA_F$
        \item $NW_F=A_F=\emptyset$.
    \end{enumerate}
\end{example}
\begin{example}[\bf A second case with infinitely many separatrices]
    \label{ex:inf seps 2}
    Let $F$ be the flow of the $C^\omega$ vector field $\xi(x,y)=(\sin y,\cos^2 y)$, whose orbits are sketched in Fig.~\ref{fig: infinitely many separatrices}(right).
    
    Notice that $\fL_\xi H(x,y)=0$ for $H(x,y)=\sec y - x$, so that 
    $$
    \cO_F(x_0,y_0)\subset
    \{(x,y): H(x,y)=H(x_0,y_0)\}.
    $$
    Like in the previous example, the separatrices of this flow are the horizontal lines $s_k=\{y=k\pi+\pi/2\}$, $k\in\bZ$.
    The relations between the separatrices, though, are quite different:
    $$
    \dots\prec s_{-1}\prec s_0\prec s_1\prec\dots.
    $$
    Hence 
    $$
    \Lambda^1_F(p)=
    \begin{cases}
        \emptyset,&p\not\in s_n\text{ for any }n\in\bZ;\\
        s_{n+1},&p\in s_n\text{ for some }n\in\bZ\\
    \end{cases}
    $$
    and
    $$
    \Lambda^2_F(p)=
    \begin{cases}
        \emptyset,&p\not\in s_n\text{ for any }n\in\bZ;\\
        \bigcup\limits_{m=n+1}^\infty s_m,&p\in s_n\text{ for some }n\in\bZ,\\
    \end{cases}
    $$
    so that: 
    \begin{enumerate}
        \item $\Lambda^\star_F(p)=\Lambda^2_F(p)$ for each 
        $p\in\bR^2$;
        \item $rank(F)=2$;
        \item $\NW_F\subsetneq\cA_F$;
        \item $NW_F=A_F=\emptyset$.
    \end{enumerate}
\end{example}
Next example shows that wandering flows can have generalized recurrent points.
\begin{example}[\bf Two wandering flows on the cylinder]
    \label{ex: cil}
    Flows of  Example~\ref{ex:inf seps 1} and Example~\ref{ex:inf seps 2} project to wandering flows on the cylinder, seen as the quotient of the plane under the identification of all pairs $(x,y)$,$(x',y')$ with $x=x'$ and $y-y'\in2\pi\bZ$.
    We denote these two flows by $F_1$ and $F_2$.
    In both cases there are exactly two separatrices, that we denote by $s_1$ and $s_2$ using the same naming of the corresponding planar examples (see Fig.~\ref{fig: infinitely many separatrices}).
    Notice that, in the quotient, all odd-numbered separatrices are equivalent to $s_1$ and all even-numbered ones to $s_2$.
    
    In case of $F_1$, the flow brings points close to $s_1$ towards $s_2$ from both sides of $s_1$.
    Hence, 
    $$
    \Lambda^1_{F_1}(p)=
    \begin{cases}
        \emptyset,&p\not\in s_1;\\
        s_2,&p\in s_1.\\
    \end{cases}
    $$
    Consequently: 
    \begin{enumerate}
        \item $\Lambda^\star_{F_1}(p)=\Lambda^1_F(p)$ for each 
        $p\in\bR^2$;
        \item $rank(F_1)=1$;
        \item $\NW_{F_1}=\cA_{F_1}$;
        \item $NW_{F_1}=A_{F_1}=\emptyset$.
    \end{enumerate}
    
    In the second case, on the contrary, points right below $s_1$ asymptote back to $s_1$ while points right above $s_1$ asymptote to $s_2$.
    As a consequence, 
    $$
    \Lambda^1_{F_2}(p)=
    \begin{cases}
        \emptyset,&p\not\in s_1\cup s_2\\
        s_2,&p\in s_1\\
        s_1,&p\in s_2\\
    \end{cases}
    $$
    and
    $$
    \Lambda^2_{F_2}(p)=
    \begin{cases}
            \emptyset,&p\not\in s_1\cup s_2\\
            s_1\cup s_{2},&p\in s_{1}\cup s_2,\\
        \end{cases}
    $$
    so that: 
    \begin{enumerate}
        \item $\Lambda^\star_{F_2}(p)=\Lambda^2_F(p)$ for each 
        $p\in\bR^2$;
        \item $rank(F)=2$;
        \item $\NW_{F_2}\subsetneq\cA_{F_2}$;
        \item $NW_{F_2}=\emptyset$;
        \item $A_{F_2}=s_1\cup s_2$.
    \end{enumerate}   
\end{example}
\section{Regular flows on $\bR^2$}
\label{sec: results}
%
%
{\BF In this section, we restrict to the case $X=\bR^2$ and $F$ will denote a regular planar flow.}
This case enjoys the following particularly strong and crucial property: 
%
\begin{proposition}
    \label{Jordan}
    A flow $F$ on $\Rt$ is wandering 
    if and only if
    the complement of each orbit of $F$ in $\bR^2$ is a pair of non-empty connected and simply connected open sets.
\end{proposition}
\begin{proof}
    Since $F$ is wandering, no orbit of $F$ can be entirely contained in any disc of finite radius, namely 
    each orbit is unbounded for both $t\to+\infty$ and $t\to-\infty$.
    Hence, each orbit $\gamma$ is the stereographic projection of a continuous loop $\hat\gamma$ on the 2-sphere.
    By Jordan's curve theorem, $\hat\gamma$ divides the sphere into two simply connected components.
    Inversely, assume that each orbit of $F$ divides the plane into two open non-empty simply connected components.
    Then there $F$ has neither periodic orbits nor fixed points so that, by the Poincar\'e-Bendixon theorem, under $F$ each point is wandering.
\end{proof}
In other words, any orbit $\gamma$ divides its complement into two sets (``sides'') so that, if $x$ and $y$ are not on the same side, 
every path from $x$ to $y$ intersects necessarily $\gamma$. Since every two distinct orbits are disjoint, this means in particular that every orbit of $F$ besides $\gamma$ lies entirely in either one of those two sets.
\begin{definition}
    Let $\alpha,\beta$ be two distinct orbits of $F$. 
    We denote by {\BF$\Delta_\alpha(\beta)$} the component of $\Rt\setminus\alpha$ containing $\beta$ and by {\BF$\Delta^*_\alpha(\beta)$} the other one.
    Moreover, we set {\BF$\Theta_{\alpha,\beta}:=\Delta_\alpha(\beta)\cap\Delta_\beta(\alpha)$}.
    We say that $\gamma$ {\bf separates} $\alpha$ from $\beta$ if $\alpha\subset\Delta^*_\gamma(\beta)$.
\end{definition}
Notice that
$$
\Delta_\alpha(\beta)\cup \alpha\cup\Delta^*_\alpha(\beta)=\Rt 
$$
and
$$
\overline{\Delta^*_\alpha(\beta)}\cup\Theta_{\alpha,\beta}\cup\overline{\Delta^*_\beta(\alpha)}=\Rt
$$
for every distinct $\alpha$ and $\beta$ and that, in both cases, the unions are disjoint.

The reader can verify the following claim.
\begin{proposition}
    For all orbits $\alpha,\beta$, the sets $\Delta_\alpha(\beta)$, $\Delta^*_\alpha(\beta)$ and $\Theta_{\alpha,\beta}$ are saturated.
\end{proposition}
%
    In order to keep the notation light, for a saturated set $A$ we will use the same symbol to denote the set of all orbits contained in $A$ and to denote the set of all points on those orbits.
    So, given an orbit $\gamma$ in $A$, if we write that $\gamma\in A$ we are considering $A$ as a collection of orbits while if we write $\gamma\subset A$ we are considering $\gamma$ and $A$ as sets of points.
\black
\begin{definition}[Kaplan, 1940~\cite{Kap40}]
    \label{def: |a,b,c|}
    Let $\cF$ be a foliation.
    Given any three mutually distinct orbits of $\cF$ denoted by $\alpha,\beta,\gamma$, we write $\alpha|\beta|\gamma$ if $\beta$ separates $\alpha$ from $\gamma$.
    Otherwise, we write $|\alpha,\beta,\gamma|^+$ if 
    there is a closed loop $\ell$ positively oriented and intersecting each orbit at a single point; similarly, we write $|\alpha,\beta,\gamma|^-$ if the orientation of $\ell$ is negative (Theorem~22 in~\cite{Kap40} grants that these cover all possible positions of a triple of distinct orbits in a regular planar flow).
\end{definition}
\begin{example}
    Consider the flow in Figure~\ref{fig:  waz}. 
    Then the orbits $s_1,s_2,s_3$ are in the chordal relation $s_1|s_2|s_3$.
    Consider now the flow in Figure~\ref{fig: finitely many seps} and let $\gamma$ be any orbit between $s_+$ and $s_-$.
    Then those three orbits are in the chordal relation $|s_-,\gamma,s_+|^+$.
\end{example}
%



Kaplan extracted from the topology of planar foliations the following set of axioms.
%
\begin{definition}[Kaplan, 1940~\cite{Kap40}]
    A {\bf chordal system} is a set $E$ endowed with two ternary relations on triplets of mutually distinct elements $a,b,c\in E$, denoted by $a|b|c$ and $|a,b,c|^+$, and satisfying the following axioms:

    \smallskip\noindent
    {\bf Axiom 1:} Either $a|b|c$ or $b|c|a$ or $c|a|b$ or $|a,b,c|^+$ or $|a,c,b|^+$.
    
    \smallskip
    From now on, we use the notation $|a,b,c|^-$ to denote $|a,c,b|^+$ and write $|a,b,c|^\pm$ to say that either $|a,b,c|^+$ or $|a,b,c|^-$.

    \smallskip\noindent
    {\bf Axiom 2.1:} $a|b|c\iff c|b|a$.

    \smallskip\noindent
    {\bf Axiom 2.2:} $|a,b,c|^+\iff|b,c,a|^+$.

    \smallskip\noindent
    {\bf Axiom 3.1:} $|a,b,c|^+,|a,c,d|^+$ $\implies$ $|a,b,d|^+,|b,c,d|^+$.

    \smallskip\noindent
    {\bf Axiom 3.2:} $|a,b,c|^\pm,\;a|b|d\implies |a,d,c|^\pm,\;c|b|d$.

    \smallskip\noindent
    {\bf Axiom 3.3:} $a|b|c,\;b|c|d\implies a|b|d,\;a|c|d$.

    \smallskip\noindent
    {\bf Axiom 3.4:} At most two of the relations $b|a|c$, $b|a|d$ and $c|a|d$ can hold.
\end{definition}
The following proposition shows that these axioms completely characterize the space of orbits of a regular planar flow.
\begin{thmX}[Kaplan, 1940~\cite{Kap40}, Theorem 25]
    The set of the orbits of a regular flow on $\bR^2$ with the relations in Definition~\ref{def: |a,b,c|} is a chordal system.
\end{thmX}

\black

\begin{proposition}
    \label{lemma: inseparable}
    Let $s_1,s_2$ two inseparable orbits of $F$.
    Then there is no orbit $\gamma$ of $F$ 
    that separates $s_1$ from $s_2$.
\end{proposition}
%
%
\begin{proof}
    Suppose that $\gamma$ does separate $s_1$ from $s_2$.
    Then $\Delta_\gamma(s_1)$ and $\Delta^*_\gamma(s_1)$ are two disjoint saturated open neighborhoods of, respectively, $s_1$ and $s_2$, against the hypothesis that $s_1$ and $s_2$ are topologically inseparable.
\end{proof}
\begin{corollary}
    \label{cor: s_1,s_2}
    Let $s_1,s_2$ a pair of inseparable separatrices and $\gamma\in\Theta_{s_1,s_2}$.
    Then either $|s_1,\gamma,s_2|^+$ or $|s_1,\gamma,s_2|^-$.
\end{corollary}
\begin{definition}
    We say that a $C^0$ arc $\gamma:\bR\to\bR^2$ is {\bf everywhere transverse} to $F$ if $\gamma$ does not intersect twice the same orbit of $F$.
\end{definition}
\begin{thmX}[Kaplan, 1940~\cite{Kap40}, Theorem 29]
    \label{lemma: ||}
    Let two curves $\gamma_1,\gamma_2$ be such that there is an arc everywhere transverse to $F$ that intersects both $\gamma_1$ and $\gamma_2$.
    Then $\gamma\in\Theta_{\gamma_1\gamma_2}$ if and only if $\gamma_1|\gamma|\gamma_2$.
\end{thmX}
\begin{corollary}
    \label{cor: Joe}
    Let $\alpha$ be an orbit of $F$.
    Then there is a saturated neighborhood $U$ of $\alpha$ such that, for every pair of distinct orbits $\gamma_1,\gamma_2\subset U$ and $\gamma\in\Theta_{\gamma_1\gamma_2}$,
    the relation $\gamma_1|\gamma|\gamma_2$ holds.
\end{corollary}
\begin{proof}
    Since $F$ is regular, every point has a neighborhood $U$ where $F$ is topologically equivalent to a linear flow.
    We can assume without loss of generality that $U$ is saturated.
    The orbits of each linear flow $L$ are a bundle of parallel lines and so $L$ admits an arc which is everywhere transverse to $L$ and cuts every of its orbits.
    Hence, there is an arc that is everywhere transverse to $F$ and cuts all orbits in $U$.
    Hence, the claim follows immediately by Lemma~\ref{lemma: ||}.
\end{proof}
\medskip\noindent
{\BF Main result \#1: planar \vectorflow flows have rank at most 2.}
We start with a few preparatory statements.
\begin{definition}
    By {\bf chain of inseparable separatrices} of $F$ we mean a finite sequence $s_1,\dots,s_k$ of orbits of $F$ such that $s_i\prec s_{i+1}$ for $i=1,\dots,k-1$.
    We say that the chain starts at $s_1$ and ends at $s_k$.
\end{definition}
\begin{corollary}
    \label{lemma: prec}
    Let $\gamma$ be an orbit of $F$ and let $s_1,\dots,s_k$ be  a chain of inseparable separatrices of $F$ such that $s_k\subset\Delta^*_\gamma(s_1)$.
    Then $\gamma=s_i$ for some $i\in\{2,\dots,k-1\}$.
\end{corollary}
\begin{proof}
    Assume that $\gamma$ is not one of the $s_i$. 
    Then there would be some $i_0\in\{2,\dots,k-1\}$ such that $s_{i_0+1}\subset\Delta^*_\gamma(s_{i_0})$,
    contradicting Lemma~\ref{lemma: inseparable}.
\end{proof}
%
%
%
\begin{proposition}
    \label{prop: L^(1,k)}
    Let $s_1,s_2$ be two orbits of $F$ and let $p\in s_1$ and $q\in s_2$.
    The following two statements are equivalent:
    \begin{enumerate}
        \item There is a chain of inseparable separatrices of $F$ starting at $s_1$ and ending at $s_2$.
        \item $q\in\Lambda^{1,k}_F(p)$ for some integer $k\geq1$.
    \end{enumerate}
\end{proposition}
\begin{proof}
    This equivalence comes immediately form the definition of $\Lambda^{1,k}$.
\end{proof}

\black
\begin{theorem}
    \label{thm1}
    Let $F$ be a regular flow on $\bR^2$. 
    Then $\Lambda^2_F(p)=\cup_k\Lambda^{1,k}_F(p)$ and 
    $\Lambda^\star_F(p)=\Lambda^2_F(p)$ for all $p\in\bR^2$.
    In particular, rank$(F)\leq2$.
\end{theorem}
\begin{proof}
    It is enough to show that $\Lambda^3_F(p)=\Lambda^2_F(p)$ for all $p\in\bR^2$.
    This is equivalent to showing that infinite chains of inseparable separatrices do not have accumulation points.
    We proceed by contradiction. 
    So, suppose that $q\in\Lambda^2_F(p)$ is a limit of an infinite sequence of inseparable separatrices.
    Then there are sequences $q_n,p_n\in\bR^2$ and $k_n\in\bN$ such that $q_n\to q$, $p_n\to p$, $q_n=\Lambda^{1,k_n}_F(p_n)$ and $k_n\to\infty$.
    In particular, since $k_n\to\infty$, either $q_n\neq q$ or $p_n\neq p$ for almost all $n$.  
    In the argument below, we assume that $q_n\neq q$ for almost all $n$. 
    We denote by $\alpha_n$ and $\beta_n$ the orbits of $F$ passing, respectively, through $p_n$ and $q_n$ and by $\alpha$ and $\beta$ the ones passing, respectively, though $p$ and $q$.
    Hence, $\alpha_n\to\alpha$ and $\beta_n\to\beta$.
    We can assume without loss of generality that $\beta_n$ ``converges monotonically'' in the sense that $\beta_{n+1}\subset\Delta_{\beta_n}(\beta)$ for all $n$. 
    Since $F$ is regular, there is a neighborhood $U_q$ of $q$ such that the restriction of $F$ to $U_q$ is topologically equivalent to a linear flow.
    We can assume without loss of generality that $U_q$ is saturated.
    Let now $N>0$ be large enough that $\beta_N\in U_q$. 
    We can assume without loss of generality that $N$ is large enough that $\alpha_n,\alpha\subset\Delta^*_{\beta_N}(\beta)$ for all $n\geq N$.
    Hence, in particular, $\alpha_{N+1}\subset\Delta^*_{\beta_N}(\beta_{N+1})$.
    Since, by hypothesis, $q_{N+1}=\Lambda^{1,k_{N+1}}_F(p_{N+1})$, this means that there is a chain of inseparable separatrices starting at $\alpha_{N+1}$ and ending at $\beta_{N+1}$.
    By Corollary~\ref{lemma: prec}, $\beta_N$ must belong to this chain.
    In particular, there is a chain of inseparable separatrices starting at $\beta_{N}$ and ending at $\beta_{N+1}$.
    Each of these pairs, therefore, must be inside $U_q$. 
    This leads to a contradiction:
    given any pair of separatrices $s_1,s_2$ and $\gamma\subset\Theta_{s_1,s_2}$, $|s_1,\gamma,s_2|^\pm$ by Corollary~\ref{cor: s_1,s_2}; on the other side, $s_1,s_2,\gamma\subset U_q$ and so, by Corollary~\ref{cor: Joe}, $s_1|\gamma|s_2$ for all $\gamma\subset\Theta_{s_1,s_2}$.
\end{proof}

\medskip\noindent
{\BF Main result \#2: each planar regular flow has a strict Lyapunov function.}
Example~\ref{ex: cil} shows that a regular flow can have generalized recurrent points. 
Next proposition shows that, nevertheless, this is not possible for {\em planar} regular flows.
%
        
%
\begin{theorem}
    \label{thm: A_F is empty}
    Let $F$ be a regular flow on $\bR^2$.
    Then $A_F=\emptyset$, namely $F$ does not have generalized recurrent points.
\end{theorem}
\begin{proof}
    Recall that $\cA_F$ is the smallest closed quasi-order containing $\NW_F$ and that $A_F$ is the set of all  points $p$ that either are fixed or for which there is another point $q$ such that $p\Ato q$ and $q\Ato p$.
    By Theorem~\ref{thm: Down A} and Theorem~\ref{thm1}, we know that $p\Ato q$ if and only if either $q$ is in the orbit of $p$ or $p$ belongs to some separatrix $s$, $q$ to some separatrix $s'$ and there is a finite chain of inseparable separatrices $s\prec\dots\prec s'$.
    Hence, in order for $p$ to be generalized recurrent, either $p$ must be fixed or there must be a chain of (at least three) inseparable separatrices starting and ending with $s$.
    This is impossible because, 
    in a chain $s_1\prec s_2\prec\dots\prec s_k$, 
    $s_1$ and $s_k$ lie in opposite 
    sides with respect to $s_2$, so that $s_k\neq s_1$. 
\end{proof}
\begin{corollary}
    \label{cor: Lyap}
    Let $F$ be a regular flow on $\bR^2$.
    Then $F$ admits a strict $C^0$ Lyapunov function.
\end{corollary}
\begin{proof}
    This is an immediate consequence of Theorem~\ref{thm: A_F is empty} and Theorem~\ref{prop: lyap2}.
\end{proof}
The result above complements Kaplan's and Whitney's Theorem~\ref{thm: first-integral}, stating that each regular planar flows has a first-integral.
\black\medskip\noindent
{\BF Main result \#3: topological classification of regular planar flows.}
Let $F$ be a planar regular flow. 
The set of all orbits of $F$ is a foliation of the plane. 
Planar foliations were classified by W. Kaplan in~\cite{Kap40,Kap41} as a concrete case application of a more general work of Whitney on families of curves. 
In his works, Kaplan refers to foliations as {\em families of curves filling the plane} and uses an ad-hoc tool he names {\em chordal relations} to prove his results.
Here we will use some of these results to achieve a consequent classification of our flows.

%


In Theorem 27 of~\cite{Kap40}, Kaplan shows that regular families of curves filling the plane are classified by their chordal systems. 
Below, we restate his claim in terms of our main objects of study.
\begin{definition}[Kaplan, 1940~\cite{Kap40}]
    We say that two planar foliations are {\bf equivalent} if one is the image of the other under a homeomorphism of the plane into itself.
    We say that they are {\bf o-equivalent} if they are equivalent under a orientation-preserving homeomorphism. 
    Two chordal systems $E_1$ and $E_2$ are {\bf isomorphic} if there is a bijection $\phi:E_1\to E_2$ such that $|a,b,c|^+$ if and only if $|\phi(a),\phi(b),\phi(c)|^+$ and $a|b|c$ if and only if $\phi(a)|\phi(b)|\phi(c)$.
\end{definition}
%
\begin{thmX}[Kaplan, 1941]
    \label{thm: CS isomorphism}
    Let $\cF_1,\cF_2$ be two planar foliations 
    with, respectively, 
    chordal systems $E_1$ and $E_2$.
    Then $\cF_1$ and $\cF_2$ are o-equivalent 
    if and only if $E_1$ is isomorphic to $E_2$.
\end{thmX}
We will need also the following straightforward generalization of Theorem~\ref{thm: CS isomorphism}.
\begin{definition}
    We say that two planar foliations are {\bf n-equivalent} if they are equivalent under a orientation-reversing homeomorphism. 
    Two chordal systems $E_1$ and $E_2$ are {\bf anti-isomorphic} if there is a bijection $\phi:E_1\to E_2$ such that $|a,b,c|^+$ if and only if $|\phi(a),\phi(b),\phi(c)|^-$ and $a|b|c$ if and only if $\phi(a)|\phi(b)|\phi(c)$.
\end{definition}
\begin{example}
    Let $\Psi:\bR^2\to\bR^2$ be any homeomorphism of the plane that switches the orientations, for instance $\Psi(x,y)=(x,-y)$.
    Given any foliation $\cF$ of the plane with chordal system $E$, denote by $\Psi_*\cF$ the foliation of the images of the leaves of $\cF$ under $\Psi$ and denote by $\Psi_*E$ the corresponding chordal system.
    Then the map $\psi:E\to\Psi_*E$ is an anti-isomorphism, since the image under $\Psi$ of each positively oriented loop is a negatively oriented loop.
\end{example}
\begin{corollary}
    \label{cor: CS anti-isomorphism}
    Let $\cF_1,\cF_2$ be two planar foliations 
    with, respectively, 
    chordal systems $E_1$ and $E_2$.
    Then $\cF_1$ and $\cF_2$ are n-equivalent if and only if $E_1$ is anti-isomorphic to $E_2$.
\end{corollary}
\begin{proof}
    Assume that $\cF_1$ and $\cF_2$ are n-equivalent.
    Then each leaf in $\cF_2$ is the image of a leaf of $\cF_1$ under a order-reversing homeomorphism $\Phi$.
    Hence, $\Phi$ switches the orientation of each loop and so $|\alpha,\beta,\gamma|^+$ if and only if $|\Phi_*\alpha,\Phi_*\beta,\Phi_*\gamma|^-$, namely $\Phi_*$ is an anti-isomorphism.

    Assume now that there is a map $\phi:E_1\to E_2$ that is an anti-isomorphism and let $\Psi$ be an orientation-reversing homeomorphism of the plane.
    Denote by $\cF_3$ the planar foliation whose leaves are the images of the leaves of $\cF_2$ under $\Psi$ and by $E_3$ the corresponding chordal system.
    
    By construction, $\Psi$ sends 
    Then $\Psi_*:E_2\to E_3$ is an anti-isomorphism, by the first part of this proof, and so the composition $\Phi_*\circ \phi:E_1\to E_3$ is a chordal systems isomorphism.
    By Theorem~\ref{thm: CS isomorphism}, there is a orientation-preserving homeomorphism $\Psi'$ that sends the orbits of $\cF_1$ into the orbits of $\cF_3$.
    Hence, the orientation-reversing homeomorphism $(\Psi')^{-1}\circ\Psi$ sends the orbits of $\cF_1$ into the orbits of $\cF_2$, namely $\cF_1$ and $\cF_2$ are n-equivalent.
\end{proof}
Next result does not appear in Kaplan's articles nor, to our knowledge, in literature but we do not consider it a major result of our work for the following two reasons.
First, it is a natural and elementary consequence of Theorem~\ref{thm: CS isomorphism} and Corollary~\ref{cor: CS anti-isomorphism}.
Moreover, it is formulated in terms of chordal systems, while we want to formulate our main results in terms of streams.
For these reasons, we index this theorem using letters, as we do for results already available in literature.
%
\begin{thmJ*}
    \label{thm: iso or anti-iso}
    Let $\cF_1,\cF_2$ be two planar foliations with, respectively, chordal systems $E_1$ and $E_2$.    
    Then $\cF_1$ is equivalent to $\cF_2$ if and only if $E_1$ and $E_2$ are either isomorphic or anti-isomorphic.
\end{thmJ*}
%

Recall that the non-wandering relation $\NW_F$ induces the binary relation $\prec$ on the set of separatrices of $F$. 
Given the family of orbits of a flow $F$,
we show below (Proposition~\ref{lemma: chordal}) that 
the chordal relations are completely determined by this binary relation.
In turn, this will allow us to state our main result.
%

\begin{figure}
    \centering
    \begin{tabular}{cc}
    \begin{tikzpicture}[scale=0.65]

   \draw[color=blue,thick] plot [hobby] coordinates {(1,5) (2.5,1.5) (4,0)};
   \draw[color=green,thick] plot [hobby] coordinates {(-0.5,5)  (-1.2,2.5)  (-3,0)};
   \draw[color=violet,thick] plot [hobby] coordinates {(3,-0.5)  (-1.2,0)  (-2,-0.5)};

   \node[left] at (-0.5,5){$\alpha$};
   \node[right] at (4,0){$\gamma$};
   \node[left,below] at (-2,-0.5){$\beta$};

    \draw[color=yellow,thick] plot [hobby] coordinates {(-2.6,1.2)  (-1.2,0)  (-0.9,-0.4)};
    \node[left] at (-2.6,1.2){$\tau_{\alpha\beta}$};

    \draw[color=yellow,thick] plot [hobby] coordinates {(-1.1,4.5)  (-.2,4.2)  (1.6,4.5)};
    \node[right] at (1.6,4.5){$\tau_{\alpha\gamma}$};

    \draw[color=yellow,thick] plot [hobby] coordinates {(3.25,1.2)  (2,0.3)  (1.6,-0.4)};
    \node[left] at (1.6,-0.4){$\tau_{\beta\gamma}$};

    \draw[color=red,thick] plot [hobby] coordinates {(3,-0.3)  (-1.2,0.2)  (-2,-0.3)};

    \draw[color=red,thick,dashed] plot [hobby] coordinates {(2.8,0.8) 
    (1.96,1.8) (1.1,3.5) (0.4,4.3) (-0.35,3.5) (-1.3,1.8) (-2.2,0.6)};

 \end{tikzpicture}
 &
 \begin{tikzpicture}[scale=0.65]

   \draw[color=blue,thick] plot [hobby] coordinates {(1,5) (2.5,1.5) (4,0)};
   \draw[color=green,thick] plot [hobby] coordinates {(-0.5,5)  (-1.2,2.5)  (-3,0)};
   \draw[color=violet,thick] plot [hobby] coordinates {(3,-0.5)  (-1.2,0)  (-2,-0.5)};

   \node[left] at (-0.5,5){$\alpha$};
   \node[right] at (4,0){$\gamma$};
   \node[left,below] at (-2,-0.5){$\beta$};

    \draw[color=yellow,thick] plot [hobby] coordinates {(-2.6,1.2)  (-1.2,0)  (-0.9,-0.4)};
    \node[left] at (-2.6,1.2){$\tau_{\alpha\beta}$};

    \draw[color=yellow,thick] plot [hobby] coordinates {(-1.1,4.5)  (-.2,4.2)  (1.6,4.5)};
    \node[right] at (1.6,4.5){$\tau_{\alpha\gamma}$};


    \draw[color=Gold3,thick] plot [hobby] coordinates {(3,-0.3)  (-1.2,0.2)  (-2,-0.3)};
    \node[left] at (-2,-0.3){$\sigma'$};
    
    \node[left] at (1.65,-0.3){$\Delta^*_{\sigma'}(\sigma)$};

    \draw[color=Gold3,thick] plot [hobby] coordinates {(0.8,5) (2.3,1.5) (3.8,0)};
    \node[left] at (3.8,0){$\sigma$};

    \node[left] at (4,2.5){$\Delta^*_{\sigma}(\sigma')$};
    \node at (0.5,2.15){$\Theta_{\sigma,\sigma'}$};

 \end{tikzpicture}\\
 (a)&(b)\\
 \begin{tikzpicture}[scale=0.65]

   \draw[color=blue,thick] plot [hobby] coordinates {(1,5) (2.5,1.5) (4,0)};
   \draw[color=green,thick] plot [hobby] coordinates {(-0.5,5)  (-1.2,2.5)  (-3,0)};
   \draw[color=violet,thick] plot [hobby] coordinates {(3,-0.5)  (-1.2,0)  (-2,-0.5)};

   \node[left] at (-0.5,5){$\alpha$};
   \node[right] at (4,0){$\gamma$};
   \node[left,below] at (-2,-0.5){$\beta$};

    \draw[color=yellow,thick] plot [hobby] coordinates {(-2.6,1.2)  (-1.2,0)  (-0.9,-0.4)};
    \node[left] at (-2.6,1.2){$\tau_{\alpha\beta}$};

    \draw[color=yellow,thick] plot [hobby] coordinates {(-1.1,4.5)  (-.2,4.2)  (1.6,4.5)};
    \node[right] at (1.6,4.5){$\tau_{\alpha\gamma}$};


    \draw[color=Gold3,thick] plot [hobby] coordinates {(3.4,-0.4) (2.5,.15) (.4,2.2) (2.5,.7) (3.6,-0.15)};
    \node[right] at (3.4,-0.4){$\sigma'$};
    \draw[color=yellow,thick] plot [hobby] coordinates {(1.5,0)  (1.55,0.3)  (1.8,0.6)};
    \node at (1.5,-0.2){$\tau$};
    

    \draw[color=Gold3,thick] plot [hobby] coordinates {(0.8,5) (2.3,1.5) (3.8,0)};
    \node[left] at (0.8,5){$\sigma$};


 \end{tikzpicture}
 &
 \begin{tikzpicture}[scale=0.65]

   \draw[color=blue,thick] plot [hobby] coordinates {(1,5) (2.5,1.5) (4,0)};
   \draw[color=green,thick] plot [hobby] coordinates {(-0.5,5)  (-1.2,2.5)  (-3,0)};
   \draw[color=violet,thick] plot [hobby] coordinates {(3,-0.5)  (-1.2,0)  (-2,-0.5)};

   \node[left] at (-0.5,5){$\alpha$};
   \node[right] at (4,0){$\gamma$};
   \node[left,below] at (-2,-0.5){$\beta$};

    \draw[color=yellow,thick] plot [hobby] coordinates {(-2.6,1.2)  (-1.2,0)  (-0.9,-0.4)};
    \node[left] at (-2.6,1.2){$\tau_{\alpha\beta}$};



    \draw[color=Gold3,thick] plot [hobby] coordinates {(0.4,4.9) (1.9,1.3) (3.4,-0.1)};
    \node[left] at (3.9,-0.3){$\sigma'$};
    \draw[color=yellow,thick] plot [hobby] coordinates {(2.3,-0.4)  (2.4,0.1)  (2.7,0.7)};
    \node at (2.3,-0.7){$\tau$};
    

    \draw[color=Gold3,thick] plot [hobby] coordinates {(0.8,5) (2.3,1.5) (3.8,0)};
    \node at (0.8,5.2){$\sigma$};


 \end{tikzpicture}\\
 (c)&(d)\\
 \end{tabular}
    \caption{Pictures showing three orbits $\alpha,\beta,\gamma$ with $|\alpha,\beta,\gamma|^+$ under different configurations. 
    The yellow arcs meet each orbit in at most a point. 
    Each other colored line is an orbit of the flow.
    (a) Each orbit is reachable from each of the other two. 
    The picture suggests that such configuration cannot happen for a regular planar flow.
    (b,c) $\alpha$ is reachable from $\beta$ and $\gamma$ while $\beta$ is not reachable from $\gamma$.
    Under this configuration, there must exist inseparable separatrices $\sigma,\sigma'$ (painted in gold) and there are two inequivalent cases shown in panels (b) and (c).
    (d) $\gamma$ is not reachable neither from $\alpha$ nor from $\beta$.
    Under this configuration there are several cases.
    All reduce to either (b) or (c) except for the one shown in this panel.
    }
    \label{fig: triples 2}
\end{figure}
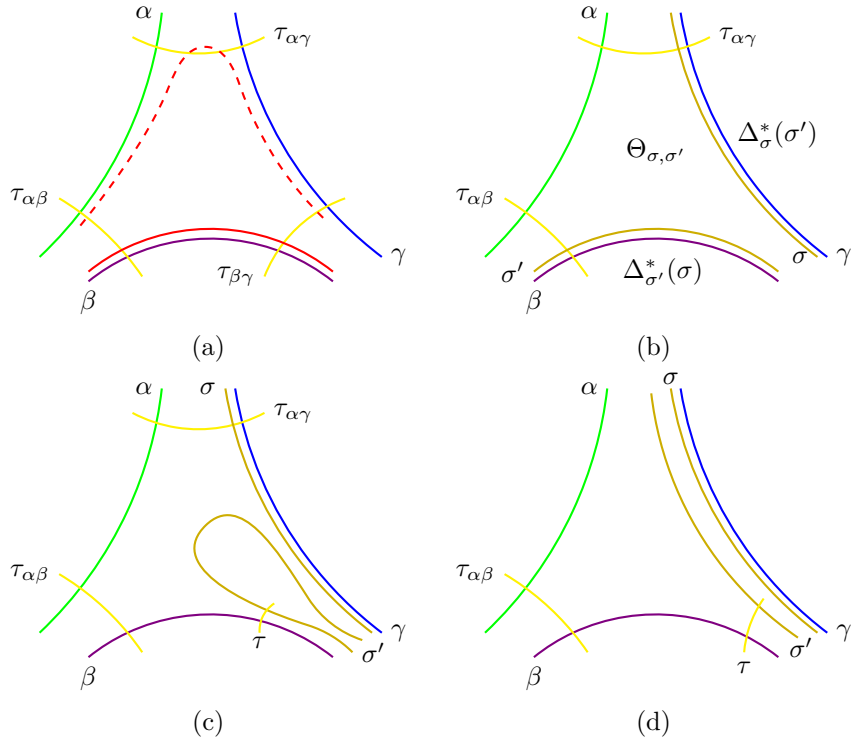

%
\begin{thmX}[Kaplan, 1940~\cite{Kap40}, Corollary to Theorem 30]
    \label{thm: arcs}
    Let $\alpha,\beta$ be two orbits of $F$ and assume that there is an arc everywhere transverse to $F$ from $p\in\alpha$ to $q\in\beta$.
    Then, there is an arc everywhere transverse to $F$ from any $p'\in\alpha$ to any $q'\in\beta$.
\end{thmX}
\begin{definition}
    Let $\alpha,\beta$ be two orbits of $F$.
    We say that {\BF $\alpha$ is reachable from $\beta$} if there is a continuous arc 
    from a point of $\alpha$ to a point of $\beta$ that meets at most once each orbit of $F$.
    We denote by $\Sigma_\alpha$ the set of all orbits of $F$ that are reachable from $\alpha$.
    We convene to include $\alpha$ in $\Sigma_\alpha$.
\end{definition}
\begin{example}
    \label{ex: sigma}
    Consider the flow whose orbits are sketched in Figure~\ref{fig: finitely many seps} (left).
    Then:
    $$
    \Sigma_\gamma=
    \begin{cases}
          \Delta_{s_2}(s_1),\ \  \gamma\in\overline{\Delta^*_{s_1}(s_2)}\\ 
          \Delta_{s_1}(s_2),\ \  \gamma\in\overline{\Delta^*_{s_2}(s_1)}\\ 
          \bR^2,\ \  \gamma\in\Theta_{s_1,s_2}\\
          
    \end{cases}
    $$
\end{example}
\begin{lemma}
    \label{lemma: reachable}
    Assume that $\alpha$ is reachable from $\beta$ and $\beta$ is reachable from $\gamma$.
    Then $\alpha$ is reachable from $\gamma$ if and only if $\beta$ separates $\alpha$ from $\gamma$, i.e. if and only if $\alpha|\beta|\gamma$.
\end{lemma}
\begin{proof}
    If $\alpha$ is reachable from $\gamma$, then $\beta\in\Theta_{\alpha\gamma}$ and so, by Theorem~\ref{lemma: ||}, $\alpha|\beta|\gamma$.
    Conversely, if $\alpha|\beta|\gamma$, by Theorem~\ref{thm: arcs} we can choose transversal arcs from $\alpha$ to $\beta$ and from $\beta$ to $\gamma$ so that the last point of the first coincides with the first point of the second.
    The new arc obtained by joining these two arcs makes $\alpha$ reachable from $\gamma$.
\end{proof}
%
In~\cite{Kap40,Kap41}, Kaplan did not consider the topology of the space of orbits; in particular, he did not consider the concept of separatrix.
Below we add a few lemmas that connect chordal relations with the prolongational relation.
\begin{lemma}
    \label{prop: Sigma}
    Let $\alpha$ be an orbit of $F$.
    The boundary of $\Sigma_\alpha$ is a cyclic collection $V_\alpha$ of separatrices of $F$.
    Moreover, for each $\sigma\in V_\alpha$, $\sigma$ is inseparable from 
    a separatrix $\sigma'\in\Sigma_\alpha\cap \Delta_\alpha(\sigma)$. 
\end{lemma}
\begin{proof}
    The fact that the boundary of $\Sigma_\alpha$ is a cyclic collection of orbits of $F$ is the content of Theorem~34 in~\cite{Kap40}.
    Here we prove that each connected component of the boundary is a separatrix which is inseparable from some separatrix inside $\Sigma_\alpha$. 
    Let $\sigma$ be an orbit at the boundary of $\Sigma_\alpha$.
    By construction, there are orbits reachable from $\alpha$ arbitrarily close to $\sigma$.
    Let $U$ be a saturated neighborhood of $\sigma$ where the flow is locally equivalent to a trivial flow and let $\gamma\in U$ be an orbit reachable from $\alpha$.
    Then there is an arc $\tau_\sigma$ from $\sigma$ to $\gamma$, since the restriction of $F$ to $U$ is trivial, and an arc $\tau_\alpha$ from $\gamma$ to $\alpha$.
    Since $\sigma\not\in\Sigma_\alpha$, these arcs cannot be combined to form an extended arc everywhere transverse to $F$.
    Hence, by Lemma~\ref{lemma: reachable}, it follows that both arcs lie in $\Delta_\gamma(\alpha)$.
    All orbits close enough to $\gamma$ in $\Delta_\gamma(\alpha)$ cut both $\tau_\sigma$ and $\tau_\alpha$.
    Let $\gamma_i$ a sequence of orbits cutting both arcs and converging to $\sigma$. 
    Then, the same sequence converges to an orbit $\sigma'$ intersecting $\tau_\alpha$, and so, in particular, reachable from $\alpha$.
    The orbits $\sigma$ and $\sigma'$ are inseparable by construction, so they are both separatrices.
    Finally, by construction, $\sigma'\in\Delta_\alpha(\sigma)$.
\end{proof}
\begin{example}
    Continuing Example~\ref{ex: sigma}, we have that
    $$
    V_\gamma=
    \begin{cases}
          s_2,\ \  \gamma\in\overline{\Delta^*_{s_1}(s_2)}\\ 
          s_1,\ \  \gamma\in\overline{\Delta^*_{s_2}(s_1)}\\ 
          \emptyset,\ \  \gamma\in\Theta_{s_1,s_2}\\
          
    \end{cases}
    $$
\end{example}
\begin{lemma}
    \label{lemma: chordal}
    Let $L$ be the space of orbits of a flow $F$ and let $\alpha,\beta,\gamma\in L$. 
    Then, modulo possibly a relabeling of the orbits,
     $|\alpha,\beta,\gamma|^+$ if and only if there are inseparable separatrices $s_1\prec s_2$ such that:
        \begin{enumerate}
            \item $\alpha\in\overline{\Delta^*_{s_2}(s_1)}$;
            \item $\beta\in\Theta_{s_1,s_2}$;
            \item $\gamma\in\overline{\Delta^*_{s_1}(s_2)}$.
        \end{enumerate}
\end{lemma}
\begin{proof}
    We prove this proposition by going over all possible combinations of which orbit is or not reachable by which orbit.

    (1) Suppose first that $\alpha$ is reachable from $\beta$ via the arc $\tau_{\alpha\beta}$, $\beta$ from $\gamma$ via $\tau_{\beta\gamma}$ and $\gamma$ from $\alpha$ via $\tau_{\gamma\alpha}$.
    Denote by $K$ the compact region bounded by $\alpha,\beta,\gamma$ and the three transversal arcs, by $K_{\alpha}$ the part of $K$ covered by orbits that intersect $\tau_{\alpha\gamma}$ and $\tau_{\alpha\beta}$ and so on.
    The sets $K_\alpha, K_\beta, K_\gamma$ are open and disjoint.
    Since $K$ is connected, there are points in $K$ that do not belong to any of them.
    Let $x$ be any of these points.
    Its orbit $\cO(x)$ cannot get out of $K$ since no orbit can cut more than twice any of the transversals.
    Hence, since $K$ is compact, $\cO(x)$ must have a non-trivial limit set. 
    This contradicts the assumption that $F$ is wandering.   

    (2) Assume now that both $\beta$ and $\gamma$ are reachable from $\alpha$ via $\tau_{\alpha\beta}$ and $\tau_{\alpha\gamma}$ respectively but $\gamma$ is not reachable from $\beta$.
    Then, by Lemma~\ref{prop: Sigma}, there exist a separatrix $\sigma\in\partial\Sigma_\beta$ such that $\gamma\in\overline{\Delta^*_\sigma(\beta)}$ and this $\sigma$ is inseparable from a separatrix $\sigma'\in\Sigma_\beta\cap\Delta_\beta(\sigma)$ (see Fig.~\ref{fig: triples 2}).
    Notice that, since $\gamma$ is reachable from $\alpha$ and $\alpha|\sigma|\gamma$, since $\alpha\in\Sigma_\beta$, then necessarily $\sigma$ must intersect $\tau_{\alpha\gamma}$.

    There are two possibilities: either $\sigma'$ intersects $\tau_{\alpha\beta}$ (Figure~\ref{fig: triples 2}(b)) or it intersects some other transversal $\tau$ (Figure~\ref{fig: triples 2}(c)).
    In the first case,  $\alpha\in\Theta_{\sigma\sigma'}$,
    $\beta\in\overline{\Delta^*_{\sigma'}(\sigma)}$ and $\gamma\in\overline{\Delta^*_{\sigma}(\sigma')}$ and we are done.
    In the second case, close enough to $\beta$ there are orbits that intersect both $\tau_{\alpha\beta}$ and $\tau$.
    Hence, there must be a separatrix $\sigma''$ inseparable from $\sigma'$ (and so from $\sigma$) that intersects both $\tau_{\alpha\beta}$ and $\tau$, so we go back to the previous case and we are done.

    (3)
    Assume now that $\beta$ is reachable from $\alpha$ via $\tau_{\alpha\beta}$ but $\gamma$ is not reachable from $\beta$ nor $\gamma$ from $\alpha$.
    Then, as in case (2), there exists $\sigma\in\partial\Sigma_\beta$ such that $\beta|\sigma|\gamma$ and a $\sigma'\in\Sigma_\beta\cap\Delta_\beta(\sigma)$ inseparable from $\sigma$.
    There are three cases for $\sigma'$.
    
    (3.1) $\sigma'$ intersects $\tau_{\alpha\beta}$.
    We are now in case (2) (see Figure~\ref{fig: triples 2}(b), except now there is no transversal $\tau_{\alpha\gamma}$) and we are done.

    (3.2) $\sigma'$ does not intersect $\tau_{\alpha\beta}$ and does not divide $\beta$ from $\gamma$.
    We are now in case (2) (see Figure~\ref{fig: triples 2}(c), except now there is no transversal $\tau_{\alpha\gamma}$) and we are done.

    (3.3) $\sigma'$ does not intersect $\tau_{\alpha\beta}$ and does divide $\beta$ from $\gamma$.
    In this case, the triple $(\alpha,\beta,\sigma')$ is in the same chordal relation as $(\alpha,\beta,\gamma)$ and is covered by case (2) above.
    

    (4)
    Assume now that $\alpha$ is not reachable from $\beta$ nor from $\gamma$ and that $\gamma$ is not reachable from $\beta$ either.
    Then there is a separatrix $\sigma\in\partial\Sigma_\beta$ such that $\gamma\in\Delta^*_\sigma(\beta)$ and that is inseparable by a separatrix $\sigma'\in\Sigma_\beta$.
    There are two main cases: either $\alpha\in\Delta^*_\sigma(\beta)$ or $\alpha\not\in\Delta^*_\sigma(\beta)$.

    In the first case, there are three subcases:

    (4.1) $\alpha|\sigma'|\beta$.
    We are now in case (2) (see Figure~\ref{fig: triples 2}(b), except now there are no transversals $\tau_{\alpha\beta}$ and $\tau_{\alpha\gamma}$) and we are done.

    (4.2) $\sigma'$ does not separate $\beta$ from $\alpha$.
    We are now in case (2) (see Figure~\ref{fig: triples 2}(c), except now there is no transversals $\tau_{\alpha\beta}$ and $\tau_{\alpha\gamma}$) and we are done.

    (4.3) $\beta|\sigma'|\gamma$.
    We are now in case (3) (see Figure~\ref{fig: triples 2}(c), except now there is no transversal $\tau_{\alpha\beta}$) and we are done.

    In the second case, the triple $(\alpha,\sigma,\gamma)$ satisfies the same chordal relation as $(\alpha,\sigma,\gamma)$, so we now repeat the discussion on this last triple. 
    There are two possibilities: either, by repeating recursively this process, we fall eventually in one of the three cases above, or we end up building a sequence $\sigma_n$ of separatrices converging to $\alpha\cup\gamma$.
    The existence of such orbits arbitrarily close to the pair of curves shows that, depending on the orientation of the flow, either $\alpha\prec\gamma$ or $\gamma\prec\alpha$.
    Moreover, by construction, $\beta\in\Theta_{\alpha\gamma}$.

    Hence, the claim holds in all possible cases.
\end{proof}
%
%
%
%
\begin{lemma}
    Let $s_1$ and $s_2$ be inseparable separatrices of a planar foliation.
    Then, given an orientation on $s_1$, the chordal relation of the ordered triple $(s_1,\alpha,s_2)$ is fully determined for every $\alpha\in\Theta_{s_1,s_2}$.
\end{lemma}
\begin{proof}
    Since $\alpha\in\Theta_{s_1,s_2}$, we know that either $|s_1,\alpha,s_2|^+$ or $|s_1,\alpha,s_2|^-$.
    Modulo homeomorphisms that preserve the orientation, there are precisely two possibilities: either the configuration of $(s_1,\alpha,s_2)$ is like in Figure~\ref{fig: determined}(left) or like in Figure~\ref{fig: determined}(right).
    The orientation of $s_1$ is opposite in the two cases, so only one is compatible with the orientation of $s_1$.
\end{proof}
\begin{definition}
    Let $\cF$ be an oriented planar foliation.
    We denote by $\Delta^+_{s_1}$ the component of $\cR^2\setminus s_1$ that lies at the right of $s_1$ with respect to the orientation on $s_1$ and by $\Delta^-_{s_1}$ the other component.
\end{definition}
Notice that the definition makes sense because {\em being at the right of $s_1$} is a condition that can be tested in any trivialization of the foliation about any point of $s_1$. 
This last lemma shows that, once an orientation is chosen, the chordal relations are dictated by the prolongational relation.
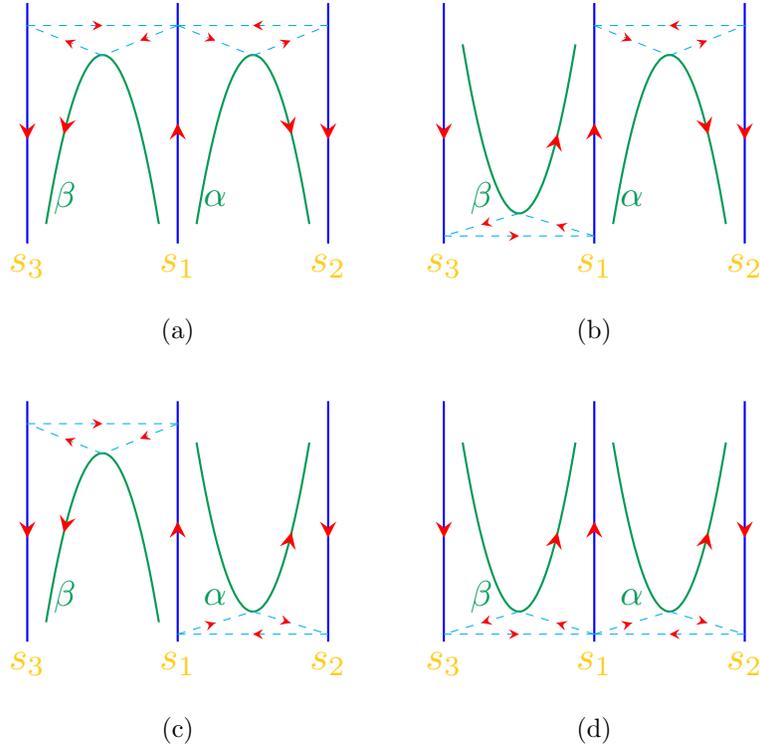
\begin{figure}
    \centering
    \tikzset{
        set arrow inside/.code={\pgfqkeys{/tikz/arrow inside}{#1}},
        set arrow inside={end/.initial=>, opt/.initial=},
        /pgf/decoration/Mark/.style={
            mark/.expanded=at position #1 with
            {
            \noexpand\arrow[\pgfkeysvalueof{/tikz/arrow inside/opt}]{\pgfkeysvalueof{/tikz/arrow inside/end}}
            }
        },
        arrow inside/.style 2 args={
            set arrow inside={#1},
            postaction={
                decorate,decoration={
                    markings,Mark/.list={#2}
                }   
            }
        },
        -<-/.style={decoration={markings,mark=at position 0.5 with {\arrow[red,>=stealth,scale=1.5]{<}}}, postaction={decorate}},
        ->-/.style={decoration={markings,mark=at position 0.5 with {\arrow[red,>=stealth,scale=1.5]{>}}}, postaction={decorate}},
        -<--/.style={decoration={markings,mark=at position 0.5 with {\arrow[red,>=stealth]{>}}}, postaction={decorate}},
    }

\definecolor{peg}{rgb}{1, 0.8, 0.1}
\setlength{\tabcolsep}{10pt} 
\renewcommand{\arraystretch}{1.5} 
\begin{tabular}{cc}

  \begin{tikzpicture}[every node/.style={scale=1.5}]
    \pgfmathsetmacro{\goldenRatio}{(1+sqrt(5)) / 2}

    \draw [color=blue,-<-,thick] (0,-0.4) -- (0,2.8);
    \draw [color=blue,->-,thick] (2,-0.4) -- (2,2.8);
    \draw [color=blue,-<-,thick] (4,-0.4) -- (4,2.8);
    \node  at (0,-0.7) {\textcolor{peg}{$s_3$}};
    \node  at (2,-0.7) {\textcolor{peg}{$s_1$}};
    \node  at (4,-0.7) {\textcolor{peg}{$s_2$}};
    \node[scale=0.9]  at (0.5,0.2) {\textcolor{ForestGreen}{$\beta$}};
    \node[scale=0.9]  at (2.5,0.2) {\textcolor{ForestGreen}{$\alpha$}};

    %
    \draw [domain=-3:3,variable=\t,smooth,samples=20,xshift=0,yshift=60,color=ForestGreen,thick]
    plot ( {-\t/4+1}, {-\t *\t/4} ) [arrow inside={end=stealth,opt={red,scale=1.5}}{0.75}];

    \draw [color=cyan,dashed] (0,2.5) -- (2,2.5)
    [arrow inside={end=stealth,opt={red,scale=1.2}}{0.5}];
    \draw [color=cyan,dashed] (1,2.1) -- (0,2.5)
    [arrow inside={end=stealth,opt={red,scale=1.2}}{0.5}];
    \draw [color=cyan,dashed] (2,2.5) -- (1,2.1)
    [arrow inside={end=stealth,opt={red,scale=1.2}}{0.5}];
    
    %
    \draw [domain=-3:3,variable=\t,smooth,samples=20,xshift=0,yshift=60,color=ForestGreen,thick]
    plot ( {\t/4+3}, {-\t*\t/4} ) [arrow inside={end=stealth,opt={red,scale=1.5}}{0.75}];
      \draw [color=cyan,dashed] (4,2.5) -- (2,2.5)
    [arrow inside={end=stealth,opt={red,scale=1.2}}{0.5}];
    \draw [color=cyan,dashed] (2,2.5) -- (3,2.1)
    [arrow inside={end=stealth,opt={red,scale=1.2}}{0.5}];
    \draw [color=cyan,dashed] (3,2.1) -- (4,2.5)
    [arrow inside={end=stealth,opt={red,scale=1.2}}{0.5}];

  \end{tikzpicture}

  &

  \begin{tikzpicture}[every node/.style={scale=1.5}]
    \pgfmathsetmacro{\goldenRatio}{(1+sqrt(5)) / 2}

    \draw [color=blue,-<-,thick] (0,-0.4) -- (0,2.8);
    \draw [color=blue,->-,thick] (2,-0.4) -- (2,2.8);
    \draw [color=blue,-<-,thick] (4,-0.4) -- (4,2.8);
    \node  at (0,-0.7) {\textcolor{peg}{$s_3$}};
    \node  at (2,-0.7) {\textcolor{peg}{$s_1$}};
    \node  at (4,-0.7) {\textcolor{peg}{$s_2$}};
    \node[scale=0.9]  at (0.5,0.2) {\textcolor{ForestGreen}{$\beta$}};
    \node[scale=0.9]  at (2.5,0.2) {\textcolor{ForestGreen}{$\alpha$}};

    %
    \draw [domain=-3:3,variable=\t,smooth,samples=20,xshift=0,yshift=0,color=ForestGreen,thick]
    plot ( {0.25*\t+1}, {0.25*\t*\t} ) [arrow inside={end=stealth,opt={red,scale=1.5}}{0.75}];

    \draw [color=cyan,dashed] (4,2.5) -- (2,2.5)
    [arrow inside={end=stealth,opt={red,scale=1.2}}{0.5}];
    \draw [color=cyan,dashed] (2,2.5) -- (3,2.1)
    [arrow inside={end=stealth,opt={red,scale=1.2}}{0.5}];
    \draw [color=cyan,dashed] (3,2.1) -- (4,2.5)
    [arrow inside={end=stealth,opt={red,scale=1.2}}{0.5}];

    %
    \draw [domain=-3:3,variable=\t,smooth,samples=20,xshift=0,yshift=60,color=ForestGreen,thick]
    plot ( {\t/4+3}, {-\t*\t/4} ) [arrow inside={end=stealth,opt={red,scale=1.5}}{0.75}];

    \draw [color=cyan,dashed] (0,-0.3) -- (2,-0.3)
    [arrow inside={end=stealth,opt={red,scale=1.2}}{0.5}];
    \draw [color=cyan,dashed] (1,0) -- (0,-0.3)
    [arrow inside={end=stealth,opt={red,scale=1.2}}{0.5}];
    \draw [color=cyan,dashed] (2,-0.3) -- (1,0)
    [arrow inside={end=stealth,opt={red,scale=1.2}}{0.5}];

  \end{tikzpicture}\\

    (a)&(b)\\
    &\\

  \begin{tikzpicture}[every node/.style={scale=1.5}]
    \pgfmathsetmacro{\goldenRatio}{(1+sqrt(5)) / 2}

    \draw [color=blue,-<-,thick] (0,-0.4) -- (0,2.8);
    \draw [color=blue,->-,thick] (2,-0.4) -- (2,2.8);
    \draw [color=blue,-<-,thick] (4,-0.4) -- (4,2.8);
    \node  at (0,-0.7) {\textcolor{peg}{$s_3$}};
    \node  at (2,-0.7) {\textcolor{peg}{$s_1$}};
    \node  at (4,-0.7) {\textcolor{peg}{$s_2$}};
    \node[scale=0.9]  at (0.5,0.2) {\textcolor{ForestGreen}{$\beta$}};
    \node[scale=0.9]  at (2.5,0.2) {\textcolor{ForestGreen}{$\alpha$}};

    \draw [domain=-3:3,variable=\t,smooth,samples=20,xshift=0,yshift=60,color=ForestGreen,thick]
    plot ( {-\t/4+1}, {-\t *\t/4} ) [arrow inside={end=stealth,opt={red,scale=1.5}}{0.75}];

    \draw [color=cyan,dashed] (0,2.5) -- (2,2.5)
    [arrow inside={end=stealth,opt={red,scale=1.2}}{0.5}];
    \draw [color=cyan,dashed] (1,2.1) -- (0,2.5)
    [arrow inside={end=stealth,opt={red,scale=1.2}}{0.5}];
    \draw [color=cyan,dashed] (2,2.5) -- (1,2.1)
    [arrow inside={end=stealth,opt={red,scale=1.2}}{0.5}];

    %
    \draw [domain=-3:3,variable=\t,smooth,samples=20,xshift=0,yshift=0,color=ForestGreen,thick]
    plot ( {0.25*\t+3}, {0.25*\t*\t} ) [arrow inside={end=stealth,opt={red,scale=1.5}}{0.75}];

    \draw [color=cyan,dashed] (4,-0.3) -- (2,-0.3) 
    [arrow inside={end=stealth,opt={red,scale=1.2}}{0.5}];
    \draw [color=cyan,dashed] (2,-0.3) -- (3,0)
    [arrow inside={end=stealth,opt={red,scale=1.2}}{0.5}];
    \draw [color=cyan,dashed] (3,0) -- (4,-0.3) 
    [arrow inside={end=stealth,opt={red,scale=1.2}}{0.5}];

  \end{tikzpicture}

  & 

  \begin{tikzpicture}[every node/.style={scale=1.5}]
    \pgfmathsetmacro{\goldenRatio}{(1+sqrt(5)) / 2}

    \draw [color=blue,-<-,thick] (0,-0.4) -- (0,2.8);
    \draw [color=blue,->-,thick] (2,-0.4) -- (2,2.8);
    \draw [color=blue,-<-,thick] (4,-0.4) -- (4,2.8);
    \node  at (0,-0.7) {\textcolor{peg}{$s_3$}};
    \node  at (2,-0.7) {\textcolor{peg}{$s_1$}};
    \node  at (4,-0.7) {\textcolor{peg}{$s_2$}};
    \node[scale=0.9]  at (0.5,0.2) {\textcolor{ForestGreen}{$\beta$}};
    \node[scale=0.9]  at (2.5,0.2) {\textcolor{ForestGreen}{$\alpha$}};

    %
    \draw [domain=-3:3,variable=\t,smooth,samples=20,xshift=0,yshift=0,color=ForestGreen,thick]
    plot ( {0.25*\t+1}, {0.25*\t*\t} ) [arrow inside={end=stealth,opt={red,scale=1.5}}{0.75}];
    \draw [color=cyan,dashed] (0,-0.3) -- (2,-0.3)
    [arrow inside={end=stealth,opt={red,scale=1.2}}{0.5}];
    \draw [color=cyan,dashed] (1,0) -- (0,-0.3)
    [arrow inside={end=stealth,opt={red,scale=1.2}}{0.5}];
    \draw [color=cyan,dashed] (2,-0.3) -- (1,0)
    [arrow inside={end=stealth,opt={red,scale=1.2}}{0.5}];

    \draw [domain=-3:3,variable=\t,smooth,samples=20,xshift=0,yshift=0,color=ForestGreen,thick]
    plot ( {0.25*\t+3}, {0.25*\t*\t} ) [arrow inside={end=stealth,opt={red,scale=1.5}}{0.75}];

    \draw [color=cyan,dashed] (4,-0.3) -- (2,-0.3) 
    [arrow inside={end=stealth,opt={red,scale=1.2}}{0.5}];
    \draw [color=cyan,dashed] (2,-0.3) -- (3,0)
    [arrow inside={end=stealth,opt={red,scale=1.2}}{0.5}];
    \draw [color=cyan,dashed] (3,0) -- (4,-0.3) 
    [arrow inside={end=stealth,opt={red,scale=1.2}}{0.5}];

  \end{tikzpicture}\\

    (c)&(d)\\

\end{tabular}
    \caption{The four possible combinations of chordal relations between the triples $(s_1,\alpha,s_2)$ and $(s_1,\beta,s_3)$, where $s_1$ is inseparable from both $s_2$ and $s_3$, $\alpha\in\Theta_{s_1,s_2}$ and $\beta\in\Theta_{s_1,s_3}$. The chordal relations are as follows: (a) $|s_1, \alpha, s_2|^+,|s_1, \beta, s_2|^-$; (b) $|s_1, \alpha, s_2|^+,|s_1, \beta, s_2|^+$; (c) $|s_1, \alpha, s_2|^-,|s_1, \beta, s_2|^-$; (d) $|s_1, \alpha, s_2|^-,|s_1, \beta, s_2|^+$.}
    \label{fig: UU}
\end{figure}
\begin{lemma}
    \label{lemma: prec => CS}
    Let $\cF$ be an oriented planar foliation and let $s_1,s_2,s_3\in\cF$ be such that $s_2,s_3$ are inseparable from $s_1$ and that $s_2\in\Delta^+_{s_1}$ and $s_3\in\Delta^-_{s_1}$.
    Then, for every $\alpha\in\Theta_{s_1,s_2}$:
    \begin{enumerate}
        \item $s_1\prec s_2$ if and only if $|s_1,\alpha,s_2|^+$;
        \item $s_2\prec s_1$ if and only if $|s_1,\alpha,s_2|^-$.
    \end{enumerate}
    Moreover, if $s_1\prec s_2$, then, for every  $\beta\in\Theta_{s_1,s_3}$:
    \begin{enumerate}
        \item $s_1\prec s_3$ if and only if $|s_1,\beta,s_3|^-$;
        \item $s_3\prec s_1$ if and only if $|s_1,\beta,s_3|^+$.
    \end{enumerate}
    While, if $s_2\prec s_1$, then, $\beta\in\Theta_{s_1,s_3}$:
    \begin{enumerate}
        \item $s_1\prec s_3$ if and only if $|s_1,\beta,s_3|^+$;
        \item $s_3\prec s_1$ if and only if $|s_1,\beta,s_3|^-$.
    \end{enumerate}   
    %
\end{lemma}   
\begin{proof}
    All possible configurations are schematically shown in Figure~\ref{fig: UU}.
    Here we discuss case (a), where $s_1\prec s_2$ and $s_1\prec s_3$, and all other cases can be done similarly.
    Let $x\in s_1$ and $y\in s_2$.
    By assumption, there are curves in $\Delta_{s_1}^+$ passing arbitrarily close to $x$ whose forward orbit passes arbitrarily close to $y$.
    We fix two charts $U_x$ about $x$ and $U_y$ about $y$ such that $U_x\cap\cF$ and $U_y\cap\cF$ are homeomorphic to a rectangle foliated by parallel vertical lines.
    Starting from $x$, we can trace a horizontal arc transversal to $\cF$ that reaches a point $z\in\alpha$.
    From $z$, we trace a small transversal arc going back towards $s_1$ and then we continue moving on a leaf of $\cF$.
    By assumption, this will lead us inside $U_y$.
    From within $U_y$, we can then continue the arc horizontally until we reach $s_2$.
    Finally, we repeat this step ``backward'' until we get back to $x$.
    This way, we build a loop that is transversal to $\cF$, touches exactly once $s_1$, $\alpha$ and $s_2$ and turns counterclockwise.
    Hence, $|s_1, \alpha, s_2|^+$.
    The same type of argument shows that $|s_1, \beta, s_3|^-$.
\end{proof}

We are now ready to prove our major result.
\begin{theorem}
    \label{thm3}
    Let $F_1,F_2$ be regular flows on the plane and denote by $\cF_1,\cF_2$ the corresponding spaces of orbits endowed with the topology induced by the projection.
    Then $F_1$ and $F_2$ are topologically equivalent if and only if there is a homeomorphism $\phi:\cF_1\to\cF_2$ that preserves the relation $\prec$, namely such that $\cO_F(x_1)\prec\cO_F(x_2)$ if and only if $\phi(\cO_F(x_1))\prec\phi(\cO_F(x_2))$ for every $x_1,x_2\in X$.
\end{theorem}
\begin{proof}
    Assume first that $F_1$ is topologically equivalent to $F_2$.
    Then there exists a homeomorphism $\Phi$ of the plane into itself that brings orbits of $F_1$ into orbits of $F_2$ and preserving the orbits orientations.
    Hence the map $\phi$ that sends $\cO_{F_1}(x)$ to $\cO_{F_2}(\Phi(x))$ is a well-defined, bijective and continuous map from $\cF_1$ to $\cF_2$.
    Moreover, since $\Phi$ preserves orientation, it induces a isomorphism between $\NW_{F_1}$ and $\NW_{F_2}$ and so $\phi$ preserves $\prec$.

    Assume now that there exists a homeomorphism $\phi:\cF_1\to\cF_2$ that preserves $\prec$.
    We claim that $\phi$ is either a chordal isomorphism or a chordal anti-isomorphism between $E_1$ and $E_2$.
    Let $\alpha,\beta,\gamma$ be three distinct orbits of $F_1$.
    If $\alpha|\beta|\gamma$, then $\alpha$ and $\gamma$ are on opposite sides with respect to $\beta$.
    Then, $\phi(\alpha)$ and $\phi(\gamma)$ are at opposite sides of $\phi(\beta)$ since, if it weren't like this then, by continuity, there would be more than one orbit whose image is $\phi(\beta)$, against the hypothesis that $\phi$ is a homeomorphism.
    Hence, $\phi(\alpha)|\phi(\beta)|\phi(\gamma)$.
    If $|\alpha,\beta,\gamma|^+$ 
    then, by Lemma~\ref{lemma: chordal}, there are inseparable separatrices $s_1,s_2$ with $s_1\prec s_2$ 
    such that $\alpha\in\overline{\Delta^*_{s_2}(s_1)}$, $\beta\in\Theta_{s_1,s_2}$ and $\gamma\in\overline{\Delta^*_{s_1}(s_2)}$, so that either $|\phi(\alpha),\phi(\beta),\phi(\gamma)|^+$ or $|\phi(\alpha),\phi(\beta),\phi(\gamma)|^-$.
    By Lemma~\ref{lemma: prec => CS}, a standard induction argument shows that, once the chordal relation value for a triple is established, all others are dictated by the $\prec$ relation.
    Hence, if $(s_1,\alpha,s_2)$ are in the same chordal relation as $(\phi(s_1),\phi(\beta),\phi(s_2))$, then all are and so $\phi$ is an isomorphism of chordal systems.
    If, on the contrary, $(s_1,\alpha,s_2)$ are in opposite chordal relation as $(\phi(s_1),\phi(\beta),\phi(s_2))$, then all are and so $\phi$ is an anti-isomorphism of chordal systems.
    
    By Theorem~J', 
    in either case there is a homeomorphism $\Phi$ of the plane into itself that brings orbits of $F_1$ into orbits of $F_2$ and viceversa. 
    The fact that $\phi$ preserves $\prec$ implies that $\Phi$ preserves the orientations of the orbits.
    Hence, $F_1$ and $F_2$ are topologically equivalent.
\end{proof}
Our main result can be reformulated as follows.
\begin{definition}
    Given a binary relation $A$ on $X$, namely $A\subset X\times X$, and a map $\Phi:X\to X$, we set 
    $$
    \Phi_*A = \{(\Phi(x),\Phi(y))\,|\,(x,y)\in A\}.
    $$ 
    Given two binary flows $F_1,F_2$ on $X$ and a continuous map $\Phi$, we say that {\BF $\Phi$ preserves the prolongational relation} (resp., the Auslander relation) if $\Phi_*\NW_{F_1}=\NW_{F_2}$ (resp., $\Phi_*\cA_{F_1}=\cA_{F_2}$).
\end{definition}
\noindent
{\bf Theorem 3'.}
{\em A regular planar flow $F_1$ is topologically equivalent to a planar flow $F_2$ if and only if there is a homeomorphims of the plane into itself that brings orbits of $F_1$ in orbits of $F_2$ and preserves the prolongation relation (equivalently, the Auslander stream).
Equivalently, each regular planar stream is completely characterized, modulo topological equivalence, by its following two topological invariants: the topology of its space of orbits and its prolongational relation (equivalently, its Auslander stream).}

\smallskip
As a final remark, we point out that the second author, jointly with Jim Yorke, showed in~\cite{DLY24} that, when the Auslander stream has countably many nodes, it coincides with the Consley's chains stream. 
Hence, all results of the present article hold after replacing Auslander's stream with the chains stream.
\section*{Acknowledgments}
The authors are thankful to Jim Yorke for several fruitful discussions about some parts of the article.
\bibliographystyle{amsplain}  
\bibliography{refs}  

\providecommand{\bysame}{\leavevmode\hbox to3em{\hrulefill}\thinspace}
\providecommand{\MR}{\relax\ifhmode\unskip\space\fi MR }
\providecommand{\MRhref}[2]{%
  \href{http://www.ams.org/mathscinet-getitem?mr=#1}{#2}
}
\providecommand{\href}[2]{#2}
\begin{thebibliography}{10}

\bibitem{Arn92b}
V.I. Arnold, \emph{Ordinary {D}ifferential {E}quations}, Springer Science \&
  Business Media, 1992.

\bibitem{Aus63}
J.~Auslander, \emph{Generalized recurrence in dynamical systems}, Contr. Diff.
  Eqs. \textbf{3} (1963), 65--74.

\bibitem{AG97}
J.~Auslander and M.~Guerin, \emph{Regional proximality and the prolongation},
  (1997).

\bibitem{AS64}
J.~Auslander and P.~Seibert, \emph{Prolongations and stability in dynamical
  systems}, Annales de l'Institut Fourier, vol.~14, 1964, pp.~237--267.

\bibitem{DeL14c}
R.~De~Leo, \emph{Weak solutions of the cohomological equation in the plane for
  regular vector fields}, Mathematical Physics, Analysis and Geometry (2014).

\bibitem{DLY24}
R.~De~Leo and J.A. Yorke, \emph{Streams and graphs of dynamical systems}, arXiv
  preprint arXiv:2401.12327 (2024).

\bibitem{HR57}
A.~Haefliger and G.~Reeb, \emph{Variet\'es (non separ\'ees) a une dimension et
  structures feuillet\'ees du plan}, Enseignement Math. \textbf{3} (1957),
  107--125.

\bibitem{Kap40}
W.~Kaplan, \emph{Regular curve-families filling the plane}, Duke Math. J.
  \textbf{7:1} (1940), 154--185.

\bibitem{Kap41}
\bysame, \emph{Regular curve-families filling the plane}, Duke Math. J.
  \textbf{8:1} (1941), 11--46.

\bibitem{Kap48}
\bysame, \emph{Topology of level curves of harmonic functions}, Trans. Am.
  Math. Soc. \textbf{63:3} (1948), 514--522.

\bibitem{Ura53}
T.~Ura, \emph{Sur les courbes d{\'e}finies par les {\'e}quations
  diff{\'e}rentielles dans l’espace {\`a} $ m $ dimensions}, Annales
  scientifiques de l'{\'E}cole normale sup{\'e}rieure, vol.~70, 1953,
  pp.~287--360.

\bibitem{Waz34}
T.~Wazewski, \emph{Sur un probleme de caractere integral relatif a l'equation
  $dz/dx +q(x, y)dz/dy = 0$}, Mathematica Cluj \textbf{8} (1934), 103--116.

\bibitem{Whi33}
H.~Whitney, \emph{Regular families of curves}, Annals of Math. \textbf{34:2}
  (1933), 244--270.

\end{thebibliography}
\end{document}